\documentclass[10pt]{amsart}

\usepackage{tgtermes}

\usepackage{a4wide}
\usepackage{bbm}
\usepackage{tikz}
\usetikzlibrary{calc}
\usetikzlibrary{arrows.meta}
\usepackage{xspace}
\usepackage{paralist}
\usepackage{wrapfig}
\usepackage{multicol}
\usepackage{quoting}

\graphicspath{{graphics/}}

\usepackage[colorlinks,linkcolor=links,citecolor=cites]{hyperref}
\definecolor{links}{rgb}{.2,.1,.5}
\definecolor{cites}{rgb}{.5,.1,.2}

\usepackage[colorinlistoftodos, bordercolor=orange, backgroundcolor=orange!20, linecolor=orange, textsize=scriptsize
]{todonotes}

\newcommand{\PP}{\mathbbm{P}}
\newcommand{\RR}{\mathbbm{R}}

\newcommand{\ZZ}{\mathbbm{Z}}

\newcommand{\TA}{\mathtt{[A]}}
\newcommand{\TB}{\mathtt{[B]}}

\DeclareMathOperator{\mGP}{GP}

\newcommand{\tdet}[1]{\mathtt{\,[#1]\,}}

\newcommand{\defn}[1]{\emph{\color{blue}#1}}

\newcommand{\bblue}[1]{\color{blue}#1\color{black}}
\newcommand{\bred}[1]{\color{red}#1\color{black}}

\DeclareMathOperator{\conv}{conv}

\DeclareMathOperator{\sign}{sign}
\DeclareMathOperator{\supp}{supp}

\DeclareMathOperator{\Gr}{Gr}

\newcommand{\bx}{{\boldsymbol{x}}}

\numberwithin{equation}{section}

\newtheorem{theorem}{Theorem}
\numberwithin{theorem}{section}
\newtheorem{lemma}[theorem]{Lemma}
\newtheorem{proposition}[theorem]{Proposition}

\theoremstyle{definition}
\newtheorem{definition}[theorem]{Definition}
\newtheorem{example}[theorem]{Example}
\newtheorem{remark}[theorem]{Remark}

\newtheorem{observation}[theorem]{Observation}

\newtheorem{convention}[theorem]{Convention}

\newcommand{\blue}[1]{{\color{blue}#1}}

\title%
[Positive Pl\"ucker tree certificates for non-realizability]%
{Positive Pl\"ucker tree certificates\\ for non-realizability}%\\ \blue{(incomplete draft)}}

\author{Julian Pfeifle}

\address{Departament de Matem\`{a}tica Aplicada, Universitat
  Polit\`{e}cnica de Catalunya}

\email{julian.pfeifle@upc.edu}

\thanks{The author was supported by the grant PID2019-106188GB-I00 from the
Spanish Ministry of Education (MEC)}

\date{December 21, 2020}
\begin{document}

\begin{abstract}
  We introduce a new method for finding a non-realizability certificate of a simplicial sphere~$\Sigma$:
  we exhibit a monomial combination of classical $3$-term Pl\"ucker relations
  that yields a sum of products of determinants that are known to be positive in any realization of~$\Sigma$;
  but their sum should vanish, contradiction.
  Using this technique, we prove for the first time
  the non-realizability of a balanced $2$-neighborly $3$-sphere constructed by Zheng,
  a family of highly neighborly centrally symmetric spheres constructed by by Novik and Zheng,
  and several combinatorial prismatoids introduced by Criado and Santos.
  The method in fact works for orientable pseudo-manifolds, not just for spheres.
\end{abstract}

\maketitle

\section{Introduction}

We start with an overview in which we assume familiarity with all the terms used,
and define them all later.
For now, recall that
a simplicial $(d-1)$-sphere~$\Sigma$, given as a list of facets, is \defn{polytopal}
if it occurs as the boundary complex of some simplicial convex $d$-polytope.
We suspect that deciding whether such a polytope actually exists is very hard,
and in some cases to be discussed later we know this to be true.
%equivalent to the ``existential theory of the reals''~\cite{richter-gebert-ziegler1995}.

We present a new technique for finding \emph{final polynomials}~\cite[Theorem~4.16]{bokowski-sturmfels-89} for deciding this question.
These certificates consist of polynomial combinations of Pl\"ucker relations on~$\Sigma$
such that the result is manifestly positive in any realization of~$\Sigma$.
The existence of such a certificate proves non-realizability, because the Pl\"ucker relations vanish on any realization of~$\Sigma$,
and so cannot combine to a positive number.

\subsection{A minimal working example}
  Sometimes just one relation suffices to prove non-realizability.

  \begin{example}
    \label{ex:minimal}
  The following list of facets defines a non-realizable $3$-sphere~$\Sigma$
  with $f$-vector $(8, 27, 38, 19)$:
  
  \begin{quote}\ttfamily\small
   +[0123] -[0124] +[0135] -[0146] +[0157] -[0167] -[0234] +[0345] -[0456] +[0567] +[1237] -[1246] -[1267] +[1357] -[2347] +[2456] -[2457] -[2567] +[3457] 
  \end{quote}

  The signs before the facets turn this list of simplices into a simplicial cycle
  representing the top homology class of this $3$-sphere. 
  Now suppose we have a realization of~$\Sigma$, and consider the 3-term Pl\"ucker relation
  \[
    0
    \ = \
    \Gamma(\mathtt{045|1267})
    \ = \
    \tdet{{0 4 5} 1 2}\tdet{{0 4 5} 6 7} - \tdet{{0 4 5} 1 6}\tdet{{0 4 5} 2 7} + \tdet{{0 4 5} 1 7}\tdet{{0 4 5} 2 6}.
  \]
  By permuting the entries inside the determinants and changing the sign accordingly, we obtain
  \begin{equation}
    \label{eq:cert-easy}
      (-1)\tdet{0 1 4 2 5}\cdot (-1)\tdet{0 5 6 7 4}
      -   \tdet{0 4 6 5 1}\cdot (-1)\tdet{2 4 7 5 0}
      +   \tdet{0 1 5 7 4}\cdot     \tdet{2 4 5 6 0}  = 0,
  \end{equation}
  where now all determinants are positive.
  
  For example, $\tdet{01425}>0$ because $\tdet{01425}=-\tdet{01245}$,
  and $\tdet{01245}$ is the ``signed slack'' of $x_5$ with respect to the facet $\tdet{0124}$
  in the supposed convex realization of~$\Sigma$;
  but the orientation of~$\mathtt{0124}$ in~$\Gamma$ is~negative by the above list.
  The other determinants can be similarly checked to be positive.

  We have obtained a sum of three positive numbers that sums to zero,
  contradicting the realizability of~$\Sigma$.
  Using the notation of Definition~\ref{def:known}
  and Convention~\ref{conv:bar} introduced later on,
  the certificate~\eqref{eq:cert-easy} reads
  \[
    \Gamma(\mathtt{045|1267})
    \ = \
    \tdet{0 1 4 2| 5} \tdet{0 5 6 7| 4} +
    \tdet{0 4 6 5| 1} \tdet{2 4 7 5| 0} +
    \tdet{0 1 5 7| 4} \tdet{2 4 5 6| 0}
    \ = \
    0.
  \]
  \hfill$\qed$
\end{example}

\subsection{Jockusch's centrally symmetric $3$-sphere on 12 vertices}

Novik and Zheng~\cite{novik-zheng20,novik2020new} recently constructed several families of centrally symmetric simplicial $3$-spheres.
\cite[Problem~5.1]{novik-zheng20} asks for properties of these spheres that pertain to realizability.
As a first application of our method, we prove the non-polytopality of the simplicial $3$-sphere~$\Delta^3_6$ with $12$~vertices and $48$~facets
which underlies their construction.
This sphere was first constructed by Jockusch~\cite{jockusch95} as part of an infinite family $\{\Delta^3_n:n\ge4\}$.

\begin{proposition}
  \label{prop:jockusch}
  Jockusch's $3$-sphere $\Delta^3_6$ is not polytopal.
\end{proposition}

\begin{proof}
An orientation of this $3$-sphere consists of the following oriented facets,
  where $\bar{x}=-x$:
  
  \begin{quote} \ttfamily\small
 $+[ 1 2 5 6 ]$ $+[ \bar{1} \bar{2} \bar{5} \bar{6} ]$ $+[ \bar{1} \bar{2} 5 6 ]$ $+[ 1 2 \bar{5} \bar{6} ]$ $+[ 2 3 5 6 ]$ $+[ \bar{2} \bar{3} \bar{5} \bar{6} ]$ $+[ \bar{2} \bar{3} 5 6 ]$ $+[ 2 3 \bar{5} \bar{6} ]$ $+[ 3 4 5 6 ]$ $+[ \bar{3} \bar{4} \bar{5} \bar{6} ]$ $+[ \bar{3} \bar{4} 5 6 ]$ $+[ 3 4 \bar{5} \bar{6} ]$ $-[ 1 \bar{4} 5 6 ]$ $-[ \bar{1} 4 \bar{5} \bar{6} ]$ $+[ 1 \bar{4} \bar{5} 6 ]$ $+[ \bar{1} 4 5 \bar{6} ]$ $-[ 1 \bar{4} \bar{5} \bar{6} ]$ $-[ \bar{1} 4 5 6 ]$ $+[ 1 2 3 5 ]$ $+[ \bar{1} \bar{2} \bar{3} \bar{5} ]$ $+[ \bar{1} \bar{2} 3 5 ]$ $+[ 1 2 \bar{3} \bar{5} ]$ $-[ 1 \bar{2} 3 5 ]$ $-[ \bar{1} 2 \bar{3} \bar{5} ]$ $-[ 1 2 4 6 ]$ $-[ \bar{1} \bar{2} \bar{4} \bar{6} ]$ $-[ \bar{1} \bar{2} 4 6 ]$ $-[ 1 2 \bar{4} \bar{6} ]$ $-[ 2 3 4 6 ]$ $-[ \bar{2} \bar{3} \bar{4} \bar{6} ]$ $-[ \bar{2} \bar{3} 4 6 ]$ $-[ 2 3 \bar{4} \bar{6} ]$ $+[ 1 \bar{3} 4 6 ]$ $+[ \bar{1} 3 \bar{4} \bar{6} ]$ $+[ 2 3 4 \bar{5} ]$ $+[ \bar{2} \bar{3} \bar{4} 5 ]$ $+[ \bar{1} 2 4 \bar{5} ]$ $+[ 1 \bar{2} \bar{4} 5 ]$ $-[ 3 4 5 \bar{6} ]$ $-[ \bar{3} \bar{4} \bar{5} 6 ]$ $-[ \bar{1} 3 5 \bar{6} ]$ $-[ 1 \bar{3} \bar{5} 6 ]$ $-[ 1 2 \bar{3} 4 ]$ $-[ \bar{1} \bar{2} 3 \bar{4} ]$ $-[ 1 2 3 \bar{4} ]$ $-[ \bar{1} \bar{2} \bar{3} 4 ]$ $+[ 1 \bar{2} 3 \bar{4} ]$ $+[ \bar{1} 2 \bar{3} 4 ]$ 
  \end{quote}

Now one single Pl\"ucker relation no longer suffices.
The non-realizability follows from the existence of the following \emph{Pl\"ucker tree},
which encodes a certain monomial combination of Pl\"ucker relations:

% [ 1 2 \bar{6} \bar{4} \bar{2} ] ([ \bar{1} \bar{2} \bar{6} \bar{4} 5 ] ([ 1 2 6 4 \bar{5} ] (GP( 1 \bar{1} \bar{6} | 2 4 \bar{4} \bar{5} )) + [ 1 \bar{4} \bar{6} \bar{5} \bar{1} ] ((-GP( 1 2 4 | \bar{1} \bar{5} 6 \bar{6} )))) + [ 1 2 \bar{5} \bar{6} \bar{1} ]*[ 1 2 6 4 \bar{5} ] ((-GP( \bar{1} \bar{4} \bar{6} | 1 \bar{2} 4 5 )))) + [ 1 2 \bar{5} \bar{6} \bar{1} ]*[ 1 2 6 4 \bar{5} ]*[ \bar{1} \bar{2} \bar{6} \bar{4} 4 ] ((-GP( 1 \bar{4} \bar{6} | \bar{1} 2 \bar{2} 5 )))
% =
% [ 1 2 \bar{6} \bar{4} \bar{2} ] ([ \bar{1} \bar{2} \bar{6} \bar{4} 5 ] ([ 1 2 6 4 \bar{5} ] (+[ 1 \bar{1} 2 4 \bar{6} ]^?[ 1 \bar{4} \bar{6} \bar{5} \bar{1} ]+[ 1 2 \bar{6} \bar{4} \bar{1} ][ \bar{1} 4 \bar{6} \bar{5} 1 ]-[ 1 2 \bar{5} \bar{6} \bar{1} ][ 1 \bar{1} 4 \bar{4} \bar{6} ]^?) + [ 1 \bar{4} \bar{6} \bar{5} \bar{1} ] (+[ \bar{1} 2 4 \bar{5} 1 ][ 1 2 6 4 \bar{6} ]+[ 1 2 6 4 \bar{1} ][ 1 2 \bar{5} \bar{6} 4 ]-[ 1 \bar{1} 2 4 \bar{6} ]^?[ 1 2 6 4 \bar{5} ])) + [ 1 2 \bar{5} \bar{6} \bar{1} ]*[ 1 2 6 4 \bar{5} ] (+[ \bar{1} \bar{2} \bar{6} \bar{4} 1 ][ \bar{1} 4 5 \bar{6} \bar{4} ]+[ 1 \bar{1} 4 \bar{4} \bar{6} ]^?[ \bar{1} \bar{2} \bar{6} \bar{4} 5 ]+[ 1 \bar{1} \bar{4} 5 \bar{6} ]^?[ \bar{1} \bar{2} \bar{6} \bar{4} 4 ])) + [ 1 2 \bar{5} \bar{6} \bar{1} ]*[ 1 2 6 4 \bar{5} ]*[ \bar{1} \bar{2} \bar{6} \bar{4} 4 ] (+[ 1 2 \bar{6} \bar{4} \bar{1} ][ 1 \bar{2} \bar{4} 5 \bar{6} ]+[ \bar{1} \bar{2} \bar{6} \bar{4} 1 ][ 1 2 \bar{6} \bar{4} 5 ]-[ 1 \bar{1} \bar{4} 5 \bar{6} ]^?[ 1 2 \bar{6} \bar{4} \bar{2} ])

\ \hfill
  \begin{tikzpicture}
    \begin{scope}[
      every node/.style={blue,rectangle,draw,rounded corners=.8ex}
      ]
      \node (G1) at (0,0) { \small$\mathtt{\Gamma(124|\bar{1}\bar{5}6\bar{6})}$ };
      \node (G2) at (4.2,0) { \small$\mathtt{\Gamma(1\bar{1}\bar{6}|24\bar{4}\bar{5})}$ };
      \node (G3) at (8.4,0) { \small$\mathtt{\Gamma(\bar1\bar4\bar6|1\bar245)}$ };
      \node (G4) at (12.8,0) { \small$-\mathtt{\Gamma(1\bar4\bar6|\bar12\bar25)}$ };

    \end{scope}
    \begin{scope}[
            every node/.style={red, circle, fill=white},
            every path/.style={very thick}
      ]
      \path (G1) edge node{\small$\mathtt{[1\bar{1}24\bar{6}]^?}$} (G2);
      \path (G2) edge node{\small$\mathtt{[1\bar14\bar4\bar6]^?}$} (G3);
      \path (G3) edge node{\small$\mathtt{[1\bar1\bar45\bar6]^?}$} (G4);
    \end{scope}
  \end{tikzpicture}    
  \hfill\

  \vspace{-3ex}
  \normalsize
  In this tree, which in this case is just a path, each edge is labeled by a determinant that occurs in both incident Pl\"ucker polynomials.
  In this example, it is also the case that the sign of these edge determinants (in red) is not known;
  but that doesn't matter, as they will be eliminated anyway.
  By successively eliminating the monomials containing these determinants 
  (in any order, by Lemma~\ref{lem:eval-order}),
  we arrive at a positive combination of monomials
  made from determinants that can be checked to be positive by the orientation above:  
\small
  \begin{align*}
    \mathtt{
    [ 1 2 \bar{6} \bar{4} \bar{2} ]
    \Big([ \bar{1} \bar{2} \bar{6} \bar{4} 5 ]
    \big([ 1 2 6 4 \bar{5} ]
    \phantom{-}
    \bblue{\Gamma( 1 \bar{1} \bar{6} | 2 4 \bar{4} \bar{5} )}\quad
    } &
    \\
    \mathtt{
    {} + [ 1 \bar{4} \bar{6} \bar{5} \bar{1} ] (-\bblue{\Gamma( 1 2 4 | \bar{1} \bar{5} 6 \bar{6} )})\big)
    } &
    \\
    \mathtt{
    {} + [ 1 2 \bar{5} \bar{6} \bar{1} ]\,[ 1 2 6 4 \bar{5} ]
    (-\bblue{\Gamma( \bar{1} \bar{4} \bar{6} | 1 \bar{2} 4 5 )})\Big)
    } &
    \\
    \mathtt{
    {} + [ 1 2 \bar{5} \bar{6} \bar{1} ]\,[ 1 2 6 4 \bar{5} ]\,[ \bar{1} \bar{2} \bar{6} \bar{4} 4 ]
    (-\bblue{\Gamma( 1 \bar{4} \bar{6} | \bar{1} 2 \bar{2} 5 )})\;\;
    }
    \\[1.5ex] =  \quad
    \mathtt{
    [ 1 2 \bar{6} \bar{4} \bar{2} ] \Big([ \bar{1} \bar{2} \bar{6} \bar{4} 5 ] \big([ 1 2 6 4 \bar{5} ]
    (\bred{[ 1 \bar{1} 2 4 \bar{6} ]^?}\,[ 1 \bar{4} \bar{6} \bar{5} \bar{1} ]+
    [ 1 2 \bar{6} \bar{4} \bar{1} ]\,[ \bar{1} 4 \bar{6} \bar{5} 1 ]-
    [ 1 2 \bar{5} \bar{6} \bar{1} ]\,\bred{[ 1 \bar{1} 4 \bar{4} \bar{6} ]^?})
    } &
    \\
    \mathtt{
    {} + [ 1 \bar{4} \bar{6} \bar{5} \bar{1} ] ([ \bar{1} 2 4 \bar{5} 1 ]\,[ 1 2 6 4 \bar{6} ]+[ 1 2 6 4 \bar{1} ]\,[ 1 2 \bar{5} \bar{6} 4 ]
    -\bred{[ 1 \bar{1} 2 4 \bar{6} ]^?}\,[ 1 2 6 4 \bar{5} ])\big)
    } &
    \\
    \mathtt{
    {} + [ 1 2 \bar{5} \bar{6} \bar{1} ]\,[ 1 2 6 4 \bar{5} ] ([ \bar{1} \bar{2} \bar{6} \bar{4} 1 ]\,[ \bar{1} 4 5 \bar{6} \bar{4} ]+
    \bred{[ 1 \bar{1} 4 \bar{4} \bar{6} ]^?}\,[ \bar{1} \bar{2} \bar{6} \bar{4} 5 ]+
    \bred{[ 1 \bar{1} \bar{4} 5 \bar{6} ]^?}\,[ \bar{1} \bar{2} \bar{6} \bar{4} 4 ])\Big)
    } &
    \\
    \mathtt{
    {} + [ 1 2 \bar{5} \bar{6} \bar{1} ]\,[ 1 2 6 4 \bar{5} ]\,[ \bar{1} \bar{2} \bar{6} \bar{4} 4 ]
    \big([ 1 2 \bar{6} \bar{4} \bar{1} ]\,[ 1 \bar{2} \bar{4} 5 \bar{6} ]+[ \bar{1} \bar{2} \bar{6} \bar{4} 1 ]\,[ 1 2 \bar{6} \bar{4} 5 ]-
    \bred{[ 1 \bar{1} \bar{4} 5 \bar{6} ]^?}[ 1 2 \bar{6} \bar{4} \bar{2} ]\big)
    }
    \\[1.5ex] =  \quad
    \mathtt{
    [ \bar{1} 2 4 \bar{5} 1 ]\,[ 1 \bar{4} \bar{6} \bar{5} \bar{1} ]\,[ 1 2 \bar{6} \bar{4} \bar{2} ]\,[ 1 2 6 4 \bar{6} ]\,[ \bar{1} \bar{2} \bar{6} \bar{4} 5 ] +
    [ 1 2 6 4 \bar{1} ]\,[ 1 \bar{4} \bar{6} \bar{5} \bar{1} ]\,[ 1 2 \bar{6} \bar{4} \bar{2} ]\,[ 1 2 \bar{5} \bar{6} 4 ]\,[ \bar{1} \bar{2} \bar{6} \bar{4} 5 ] 
    } &
    \\
    \mathtt{ {}+
    [ 1 2 \bar{6} \bar{4} \bar{1} ]\,[ 1 2 \bar{5} \bar{6} \bar{1} ]\,[ 1 2 6 4 \bar{5} ]\,[ 1 \bar{2} \bar{4} 5 \bar{6} ]\,[ \bar{1} \bar{2} \bar{6} \bar{4} 4 ] +
    [ 1 2 \bar{6} \bar{4} \bar{1} ]\,[ \bar{1} 4 \bar{6} \bar{5} 1 ]\,[ 1 2 \bar{6} \bar{4} \bar{2} ]\,[ 1 2 6 4 \bar{5} ]\,[ \bar{1} \bar{2} \bar{6} \bar{4} 5 ] 
    } &
    \\
    \mathtt{ {}+
    [ 1 2 \bar{5} \bar{6} \bar{1} ]\,[ \bar{1} \bar{2} \bar{6} \bar{4} 1 ]\,[ 1 2 \bar{6} \bar{4} \bar{2} ]\,[ 1 2 6 4 \bar{5} ]\,[ \bar{1} 4 5 \bar{6} \bar{4} ] +
    [ 1 2 \bar{5} \bar{6} \bar{1} ]\,[ \bar{1} \bar{2} \bar{6} \bar{4} 1 ]\,[ 1 2 6 4 \bar{5} ]\,[ 1 2 \bar{6} \bar{4} 5 ]\,[ \bar{1} \bar{2} \bar{6} \bar{4} 4 ]
    \rlap{\ .}
      } &
  \end{align*}

  \normalsize
The existence of this combination proves the non-realizability of Jockusch's sphere~$\Delta^3_6$,
because in any polytopal realization, all Pl\"ucker relations vanish,
but all monomials in the final combination are positive.
Moreover, the method detailed in Section~\ref{sec:ip} shows that this certificate is as short as possible.
\end{proof}

Below, we build on this result and % in  below,
show that \emph{none} of Jockusch's spheres~$\Delta^3_n$ (Theorem~\ref{thm:jockusch}), and indeed
\emph{none} of the spheres in Novik \& Zheng's family $\{\Delta^d_n: n-2\ge  d\ge3\}$ (Theorem~\ref{thm:nz4})  are realizable.
For now,
after having seen the method at work in some examples,
let's start again and carefully define all terms.

\section{Motivation and overview}

Why do we even care about realizability?

\subsection{There are so many spheres. Which ones are actually polytopes?}

Write $s_d(n)$~for the number of combinatorial types of $d$-dimensional simplicial spheres on $n$~vertices,
and $p_d(n)$~for the number of combinatorial types of $d$-dimensional simplicial polytopes on $n$~vertices.
Goodman and Pollack~\cite{goodman-pollack1986},~\cite{alon1986} showed, suprisingly, that asymptotically there are only ``very few'' types of simplicial polytopes:
in fixed $d$, 
\[
  p_d(n)
  \ \in\
  2^{\Theta_d(n\log n)}.
\]
Kalai~\cite{kalai1988} used Stanley's Upper Bound Theorem~\cite{stanley1975} for simplicial spheres and his own construction of ``squeezed spheres''
to establish astronomically larger upper and lower bounds for simplicial spheres:
\[
  2^{\Omega(n^{\lfloor d/2\rfloor})}
  \ \le \
  s_d(n)
  \ \le \
  2^{O(n^{\lceil d/2\rceil}\log n)}.
\]
For even dimension~$d$, these bounds already coincide up to a term of $\log n$ in the exponent.
For $d=3$, the lower bound was first improved from~$2^{\Omega(n)}$ to~$2^{\Omega(n^{5/4})}$ by Pfeifle and Ziegler~\cite{pfeifle-ziegler2004}.
In general odd dimension $d=2k-1$, the best result known to date is due to Nevo, Santos and Wilson~\cite{nevo-santos-wilson2016}, who improved the lower bound from~$2^{\Omega(n^{k-1})}$ to $2^{\Omega(n^k)}$, asymptotically very close to the upper bound apart from, again, a factor of~$\log n$ in the exponent.

\smallskip
As stated in the introduction, we suspect that deciding whether a simplicial complex~$\Sigma$ that is guaranteed to be homeomorphic to a $(d-1)$-sphere
is actually the boundary of a convex $d$-polytope is a very hard problem.
Already for $d=4$ (in the non-simplicial case), it is equivalent to the ``existential theory of the reals''~\cite{richter-gebert-ziegler1995},
which is known to be NP-hard by the results of Mn\"ev~\cite{mnev1988} and Shor~\cite{shor1991}.
According to Adiprasito and Padrol~\cite{AdiprasitoPadrol2}, the realization problem for neighborly spheres is ``universal'',
i.e., for every ``primary basic open semi-algebraic set over $\ZZ$'',  there exists a neighborly polytope with that realization space.
But universality for simplicial polytopes in any fixed dimension~$d$ remains open.
Therefore, any new technique is welcome.

\subsection{Pl\"ucker relations}

A common strategy to prove non-realizability of a simplicial sphere $\Sigma$ is by
exhibiting a property that $\Sigma$~should have if it were realizable,
but showing that $\Sigma$~does not in fact have it.

One popular property of this kind is that realizable spheres should satisfy the Pl\"ucker relations.
For this, suppose that $\Sigma$ is realized as the boundary complex of a simplicial convex polytope~$P\subset\RR^d$ with $n$~vertices,
embed $P$ into the hyperplane $\{x=(x_0,x_1,\dots,x_d)\in\RR^{d+1}:x_0=1\}$,
and write down the resulting $(d+1)\times n$ real matrix~$V$ of homogeneous coordinates of~$P$.
The row space of~$V$ corresponds to a point on the Grassmannian~$\Gr(d+1,n)$.
After embedding into~$\RR\PP^{\binom{[n]}{d+1}}$,
this point has \defn{Pl\"ucker coordinates}
$
  \big(\det V_J : J\in \binom{[n]}{d+1}\big),
$
consisting of all $(d+1)\times(d+1)$-minors $V_J$ of~$V$
enumerated in some fixed way, for example lexicographically.
These coordinates satisfy the \emph{Pl\"ucker relations}~\cite[Theorem~1.8]{bokowski-sturmfels-89}, \cite[Proposition~2.2.10]{maclagan-sturmfels2015},
\cite[Theorem~14.6]{miller-sturmfels-2005}.
In particular,
any ordered $(d+3)$-tuple of indices $(S|ijkl)$ with $S\in\binom{[n]}{d-1}$ and $i,j,k,l\in[n]\setminus S$
gives rise to the \defn{3-term Pl\"ucker relation}
\[
  \Gamma(S|ijkl)
  \ = \
  0
\]
for~$\Sigma$, where
\[
  \bblue{\Gamma(S|ijkl)}
  \ := \
  [Sij] [Skl] - [Sik] [Sjl] + [Sil] [Sjk].
\]

We refer to \cite{gouveia-macchia-wiebe20} for an overview of four different models of realization spaces of polytopes,
which in particular relates this discussion to \emph{slack ideals} and the \emph{slack variety}~\cite{GMTW-2019} \cite{BW-2019}.

\subsection{Overall strategy}

Our strategy for proving the non-realizability of a simplicial complex~$\Sigma$
by finding a final polynomial is analogous to the use of the \emph{Positivstellensatz}~\cite{sturmfels-2002}, \cite[Chapter~4]{bokowski-sturmfels-89}
to prove the non-existence of a solution to a system of real polynomial equations;
see~\cite{deloera-malkin-parrilo-2012} and e.g.~\cite{dressler-iliman-dewolff2017}.

Namely, we view the 3-term Pl\"ucker relations as polynomials in the ring
$\RR[\bx]:=\RR\big[x_J:J\in\binom{[n]}{d+1}\big]$.
We then find and exhibit a polynomial combination $\tau = \sum \kappa_m \Gamma_m$ of certain $\Gamma_m=\Gamma(S|ijkl)$ 
such that
\begin{itemize}
\item each coefficient $\kappa_m\in\RR[\bx]$ is a monomial that manifestly only takes on positive values for any value of $\bx$
  that comes from a realization of~$\Sigma$;
\item and the combination $\tau\in\RR[\bx]$ itself is manifestly a positive sum of monomials
  that \emph{also} only take on positive values for any value of $\bx$ that comes from a realization of~$\Sigma$.
\end{itemize}
The contradiction to realizability then arises from the fact that each $\Gamma_m$ \emph{vanishes}
for any value of $\bx$ that comes from a realization of~$\Sigma$;
and we cannot get a positive value from a polynomial combination of zeroes.

We find these combinations by evaluating a certain tree whose nodes are Pl\"ucker polynomials;
hence the name \emph{positive Pl\"ucker tree certificates}.

\section{Technical apparatus}

\subsection{Notation}
Fix an integer $n\ge3$. % and $d$ such that $n-1\ge d\ge2$.
We set $\defn{$[n]$}=\{1,2,\dots,n\}$, and
take the \defn{size} of a set~$S$ to be is its cardinality~$|S|$.
An \defn{(abstract) simplicial complex} $\Sigma$ on $n$~vertices is
a collection of subsets $\sigma\subset[n]$,
called \defn{faces}, such that the collection is closed under taking subsets.
An inclusion-maximal face is called a \defn{facet}.
We will always assume that $\Sigma$ is \defn{pure},
which means that all facets have the same cardinality.
The \defn{dimension} of a face~$G$ of~$\Sigma$ is $\dim G=|G|-1$,
and the \defn{dimension} of~$\Sigma$ is the dimension of any facet.
We set $\defn{d}:=\dim\Sigma+1$, so that each facet of~$\Sigma$ has cardinality~$d$,
and $n\ge d+1\ge 3$.
A \defn{ridge} of~$\Sigma$ is a face of dimension~$(\dim\Sigma-1)$.
The \defn{$k$-skeleton} $\Sigma^k$ of~$\Sigma$ is the set of its $k$-dimensional faces.
The \defn{boundary}~$\partial\Sigma$ of~$\Sigma$ is the set of ridges contained in exactly one facet.

Fix a $(k-1)$-dimensional face $\sigma\subset[n]$.
Any permutation $\vec{\sigma}=[i_1,\dots, i_k]$ of the vertices of~$\sigma$
is called %an \defn{orientation} of~$\sigma$, or
an \defn{oriented simplex}.
The \defn{support} of~$\vec\sigma$ is $\defn{$\supp\vec\sigma$}=\{i_1,\dots,i_k\}$.
We will not obsess with notation, and happily confuse~$\sigma$ with~$\vec\sigma$ unless we fear misunderstanding.

A \defn{realization} of~$\Sigma$ in~$\RR^d$ is a set $X(\Sigma)=\{x_1,\dots,x_n\}$ of points in~$\RR^d$ such that
for each face~$\sigma\in\Sigma$, the set $\conv X_\sigma$ is a \defn{face of the convex hull} $\conv X$,
where $X_\sigma:=\{x_i:i\in\sigma\}$.
By this we mean the existence of an oriented hyperplane
$H_\sigma=\{x\in\RR^d: \langle a_\sigma,x\rangle = b_\sigma\}\subset\RR^d$
such that $\conv X_\sigma\subset H_\sigma$
and the remainder $\conv X_{[n]\setminus\supp\sigma}$ lies strictly to one side of~$H_\sigma$.
We say that $\Sigma$~is \defn{non-realizable in~$\RR^d$} if no such assignment of points in~$\RR^d$ exists.

Since we attempt to realize the complex~$\Sigma$ in~$\RR^d$,
we are especially interested in the case $k=d+1$,
which corresponds to full-dimensional simplices in~$\RR^d$.
For this, we agree to use homogeneous coordinates for our points, so that
$\sign\det(x_{i_1},\dots,x_{i_{d+1}})$ is the orientation of a full-dimensional ordered simplex in~$\RR^d$.

We will need both the sign of the determinant $\det\sigma=\det(x_{i_1},\dots,x_{i_{d+1}})$
determined by the permutation~$\vec{\sigma}$,
as well as the signature of the permutation~$\vec{\sigma}=[i_1,\dots,i_d]$ itself.
We therefore adopt the mnemonic that
\defn{$\nu(\sigma)$} is the sig\textbf{\color{blue}{\underline{n}}} of~$\det\sigma$,
and $\defn{$\varepsilon(\sigma)$}=(-1)^{t(\vec\sigma)}$ stands for the signatur\textbf{\color{blue}{\underline{e}}} of~$\vec\sigma$;
here \defn{$t(\vec\sigma)$} is the number of transpositions needed to build~$\vec\sigma$ from the identity permutation.
The two possible values~$\pm1$ of~$\varepsilon(\sigma)$ correspond to the two equivalence classes of orientations of~$\sigma$.

An \defn{oriented simplicial $k$-chain}~\cite{munkres84} for $0\le k\le \dim\Sigma$ is a formal sum
$\sum_{\sigma\in\Sigma}c(\sigma)\vec\sigma$ of oriented $k$-dimensional simplices, with coefficients $c(\sigma)\in\ZZ$.
We impose that $c(\vec\sigma) = -c(-\vec\sigma)$
whenever $-\vec\sigma$ is an oriented simplex of the opposite orientation as~$\sigma$.
The group of oriented simplicial $k$-chains is denoted~$C_k(\Sigma)$.
For $k<0$ or $k>\dim\Sigma$ we set $C_k(\Sigma)=0$.
The \defn{boundary~$\partial\vec\sigma$} of an oriented simplex is the image of~$\vec\sigma$
under the map $\partial_k:C_k(\Sigma)\to C_{k-1}(\Sigma)$ given by
\[
  \partial_k\vec\sigma
  \ = \
  \partial[i_0,\dots,i_k]
  \ = \
  \sum_{j=0}^k(-1)^j [i_0,\dots,\widehat{i_j},\dots,i_k].
\]
The kernel of $\partial_k:C_k(\Sigma)\to C_{k-1}(\Sigma)$
is the group~$Z_k(\Sigma)$ of~\defn{$k$-cycles},
while the image of $\partial_{k+1}:C_{k+1}(\Sigma)\to C_k(\Sigma)$
is called the group~$B_k(\Sigma)$ of~\defn{$k$-boundaries}.
A standard calculation shows that $\partial_k\circ\partial_{k+1}=0$, so that
$B_k(\Sigma)\subset Z_k(\Sigma)$, i.e., every boundary of a $(k+1)$-chain is automatically
a $k$-cycle.
The \defn{$k$-th integral homology group of~$\Sigma$} is
$H_k(\Sigma;\ZZ) := Z_k(\Sigma)/B_k(\Sigma)$.

\smallskip
A $(d-1)$-dimensional simplicial complex~$\Sigma$ is 

\begin{itemize}[\quad$\triangleright$]
\item
  a \defn{closed pseudomanifold} if each ridge is incident to exactly two facets;
\item
  \defn{orientable}
  if it has the top-dimensional integral homology group of a sphere, $H_{d-1}(\Sigma;\ZZ)\cong\ZZ$;
\item
  a \defn{homology sphere} if \emph{all} (reduced) homology groups coincide with those of a
  sphere, and the same holds true for all vertex links;
\item
  and a \defn{simplicial sphere} if it is homeomorphic to a
  sphere.
\end{itemize}
Thus, any simplicial sphere is a pure, closed, orientable
pseudomanifold and a homology sphere.

% A \defn{realization in~$\RR^d$} of~$\Sigma$ is a bunch of coordinates $x_i=(x_{i1},\dots,x_{id})\in\RR^d$
% for each $i=1,\dots,n$, such that the boundary complex of the convex hull of the~$x_i$ contains~$\Sigma$,
% and equals~$\Sigma$ if $\partial\Sigma=\emptyset$.

If $\Sigma$~is closed and orientable, an \defn{orientation} of~$\Sigma$ is a map $\omega:\Sigma^{d-1}\to\{+,-\}$
such that the simplicial chain $\sum_{\sigma\in\Sigma^{d-1}}\omega(\sigma)\vec\sigma$ generates $H_{d-1}(\Sigma;\ZZ)$.
If $\Sigma$~is not closed, let $C_1,\dots,C_k$ be the connected components of (the dual graph of)~$\partial\Sigma$,
and let $y_1,\dots,y_k$ be new vertices.
We extend~$\Sigma$ to an enlarged complex~$\widetilde\Sigma$ by coning over the boundary, i.e.,
adding all facets of the cones $y_j\star C_j$.
In the case $H_{d-1}(\widetilde\Sigma;\ZZ)=\ZZ$, we define the \defn{orientation} of~$\Sigma$ to be
the restriction of a generating cycle of~$H_{d-1}(\widetilde\Sigma;\ZZ)$ to~$\Sigma$,
and we say that $\Sigma$ is an \defn{orientable pseudomanifold with boundary}.

\begin{convention}
  We will always assume that $\Sigma$~is an orientable, but not necessarily closed, pseudomanifold of dimension $\dim\Sigma=d-1$.
\end{convention}

\subsection{Nonrealizability certificates from Pl\"ucker relations}

Recall that our goal is to exploit positivity
to prove that a $(d-1)$-dimensional simplicial complex~$\Sigma$
given to us as a list of facets cannot be realized in convex position.
Each 3-term Pl\"ucker relation
\[
  \Gamma(S|ijkl)
  \ = \
  [Sij] [Skl] - [Sik] [Sjl] + [Sil] [Sjk]
\]
is a linear combination of products of determinants.
If the signs of the determinants conspire together in such a way
that all three terms in this combination have a positive coefficient
(or all three have a negative coefficient), then the equation $\Gamma(S|ijkl)=0$ yields a contradiction to realizability:
no sum of positive numbers (or of negative numbers) can be zero.

\smallskip
How can we control the signs of these determinants?
Of course, we cannot control any of them individually,
but in a convex realization we do know that many of them are \emph{equal}.
To see this, suppose we have a convex realization $X(\Sigma)$ of~$\Sigma$,
and consider any oriented facet~$\vec F=(i_1,\dots,i_d)$ of~$\Sigma$.
Because the realization is convex, all points~$x_j$ with $j\notin\supp\vec F$ lie on the same side of the oriented hyperplane~$H_F$.
Another way of saying this is that the determinant~$\det(x_{i_1},\dots,x_{i_d},x_j)$
has the same sign for all indices $j\in[n]\setminus(\supp\vec F)$,
and this sign is determined by the orientation of~$\vec F$.

\smallskip
To turn this dream into a working algorithm, we first need a normal form for our determinants.

\begin{definition}[normal form]
  \label{def:known}
  Fix an oriented, not necessarily closed pseudomanifold~$\Sigma$ with orientation $\omega:\Sigma^{d-1}\to\{+,-\}$.
  An \defn{(oriented) solid} in~$\Sigma$ is an oriented simplex $\vec\pi=[i_1,\dots,i_{d+1}]$ of size~$d+1$.
  We say that an oriented solid~$\vec\pi$ is~\defn{known} (i.e., we can control its orientation) \defn{via~$F$}
  if there exists a facet~$F\in\Sigma^{d-1}$ such that $F\subset\supp\vec\pi$,
  where $F$~is not necessarily uniquely determined.
  In this case, the facet~$F$ is said to \defn{determine}~$\pi$,
  and the vertex of~$\vec\pi$ \defn{associated to~$F$ in~$\pi$} is $\defn{$\alpha(\pi,F)$}=(\supp\vec\pi)\setminus F$.
  For any oriented simplex~$\sigma$,  write \defn{$\sigma_<$} for the reordering of~$\sigma$ into ascending order.

  \smallskip
  The \defn{normal form} $\overline\pi$ of~$\vec\pi$ is defined as follows:
  \begin{itemize}[\quad$\triangleright$]
  \item Suppose first that $\vec\pi$ is known via the facet~$F\subset\Sigma^{d-1}$.
    Without loss of generality, suppose
    that~$\vec F=\vec F_<=(i_1,\dots,i_d)$ with $i_1<\dots<i_d$,
    and that the vertex of $\vec\pi$ associated to~$F$ is~$\alpha(\pi,F)=i_{d+1}$.
    \begin{itemize}[\quad$\bullet$]
    \item If $\omega(\vec F)=+$, then $\overline\pi=[i_1,\dots,i_{d-2},i_{d-1},i_d|i_{d+1}]$ (see Convention~\ref{conv:bar}).
    \item If $\omega(\vec F)=-$, then $\overline\pi=[i_1,\dots,i_{d-2},i_d,i_{d-1}|i_{d+1}]$.
    \end{itemize}
    In particular, $\varepsilon(\overline{\pi}) = \omega(\vec F_<)\varepsilon([F_<,\alpha(\pi,F)])$.

  \item If $\vec\pi$ is not known, then $\overline\pi$ is the lexicographically ordered version of~$\vec\pi$:
    $\overline\pi=[j_1,\dots,j_{d+1}]$, where $j_1<\dots<j_{d+1}$ and $\{j_1,\dots,j_{d+1}\}=\supp\vec\pi$.
  \end{itemize}
\end{definition}

\begin{convention}
 \label{conv:bar}
 We sometimes typographically distinguish the normal form of a known solid~$\pi$ (from that of an  ``unknown'' solid)
 by means of the vertical bar \defn{`` $|$''} separating~$F$ from~$\alpha(\pi,F)$.

 Sometimes, to emphasize that a solid $\overline\pi=[j_1,\dots,j_{d+1}]$ is not known, we add a question mark as a superscript:
  $\overline\pi=\bred{[j_1,\dots,j_{d+1}]^?}$. 
\end{convention}

\begin{proposition}
  \label{prop:normal-form}
  \begin{enumerate}[\upshape(1)]
  \item
    The signature $\varepsilon(\overline\pi)$ of the normal form of an oriented solid $\vec\pi$ is well-defined,
    and $\varepsilon(\overline{\pi}) = \omega(F_<)\varepsilon(\pi)$ if $\pi$~is known via $F$.

  \item 
    If $\vec\pi$ is known with normal form $\overline\pi = [i_1,\dots,i_d|i_{d+1}]$,
    we may assume withour loss of generality that $\nu(\overline{\pi}) = \sign\det(x_{i_1},\dots,x_{i_{d+1}})=+$ in any realization of~$\Sigma$.

  \item
    For any oriented solid, $\pi = \varepsilon(\pi)\varepsilon(\overline\pi)\overline\pi$.
    
  \end{enumerate}
\end{proposition}

Before proving this, let us interject that
if $\pi$ is known, then $\nu(\pi) = \varepsilon(\pi)\varepsilon(\overline\pi)\nu(\overline\pi) = \varepsilon(\pi)\varepsilon(\overline\pi)$
by Proposition~\ref{prop:normal-form}~(2) and~(3).
However, if $\pi$~is not known, we have no control over $\nu(\pi)$.
We formalize this state of affairs as follows.

\begin{definition}[$\chi$]
  If an oriented solid $\pi$~is known, set $\chi_{\pi}:=\chi(\pi):=\varepsilon(\pi)\varepsilon(\overline\pi)$,
  otherwise set  $\chi_{\pi}:=\,\bred{?}$.
  In all cases, we set $\chi_{\pi}^2=+$.
\end{definition}

We can therefore express the relationship between an oriented solid and its normal form as
\[
  \pi
  \ = \
  \chi_{\pi}\,
  \overline\pi.
\]

\begin{proof}[Proof of Proposition~\ref{prop:normal-form}]
  (1) For showing the well-definedness of $\varepsilon(\overline\pi)$,
  we may suppose that $\supp\pi=\{1,2,\dots,d+1\} = F\cup\{k\} = G\cup\{\ell\}$,
  and write $R=F\cap G=\{1,2,\dots,d+1\}\setminus\{k,\ell\}$ for the ridge common to $F$~and~$G$.
  The two representations of~$\pi$ are then as follows (suppose wlog $\ell < k$):
  \begin{center}
    \begin{tabular}[c]{c|cccc|cccc|ccc|c}
      \small position:   & \small 1 & \small 2 & \small \dots & \small $\ell-1$ & \small $\ell$ & \small $\ell+1$ & \small \dots & \small $k-1$ & \small $k$ & \small \dots & \small $d$ & \small $d+1$ \\\hline
      $F_< k$:    & 1 & 2 & \dots & $\ell-1$ & $\ell$ & $\ell+1$ & \dots & $k-1$ & $k+1$ & \dots & $d+1$ & $k$ \\
      $G_<\ell$:  & 1 & 2 & \dots & $\ell-1$ & $\ell+1$ & $\ell+2$ & \dots & $k$ & $k+1$ & \dots & $d+1$ & $\ell$ 
    \end{tabular}
  \end{center}
  We have written $F_<$, $G_<$ to emphasize that $F$ and $G$ are sorted.
  By definition,
  \begin{equation}
    \label{eq:can1}
    \varepsilon(\overline{F_<k})
    \ = \
    \omega(F_<)\varepsilon(F_<k)
    \qquad \text{and}\qquad
    \varepsilon(\overline{G_<\ell})
    \ = \
    \omega(G_<)\varepsilon(G_<\ell),
  \end{equation}
  and we wish to show that these quantities are equal.
  Now
  \begin{equation}
    \label{eq:can2}
    \omega(F_<)
    \ = \
    (-1)^{d-\ell}\omega(R_<\ell)
    \qquad\text{and}\qquad
    \omega(G_<)
    \ = \
    (-1)^{d-k+1}\omega(R_<k),
  \end{equation}
  because to sort $F_<$ into the form $R_<\ell$, we must swap the index $\ell$ to the right a total of $d-(\ell+1)+1$ steps.
  Similarly, to sort $G_<$ into the form $R_<k$, we must swap the index $k$ to the right a total of $d-k+1$ steps.
  Next,
  \begin{equation}
    \label{eq:can3}
    \varepsilon(F_<k)
    \ = \
    (-1)^{d-k+1}
    \qquad\text{and}\qquad
    \varepsilon(G_<\ell)
    \ = \
    (-1)^{d-\ell+1},
  \end{equation}
  because to sort $F_<k$ we must bring~$k$ to the left a total of $d-k+1$ steps,
  and to sort $G_<\ell$ we must bring~$\ell$ to the left a total of $d-\ell+1$ steps.
  Finally,
  \begin{equation}
    \label{eq:can4}
    \omega(R_<\ell)
    \ = \
    -\omega(R_<k)
  \end{equation}  
  because $\omega$ is an orientation.
  Substituing~\eqref{eq:can2}, \eqref{eq:can3}~and~\eqref{eq:can4}
  into \eqref{eq:can1} completes the proof.

  (2) Since $\overline\pi$ is sorted such that the facet $F$ comes first,
  and $F$ indexes a facet of the convex hull in any realization of~$\Sigma$,
  the sign $\nu(\overline\pi)$ is the same for all choices of $i_{d+1}\in[n]\setminus \supp F$.
  Moreover, since $\omega$~is an orientation, this sign is consistent across all facets~$F$ of~$\Sigma$.
  After a global reflection in~$\RR^d$, if necessary, we may therefore assume that $\nu(\overline\pi)=+$.

  (3) The sign $\varepsilon(\pi)$ reflects the change in signature as we permute $\pi$ to the identity,
  and $\varepsilon(\overline\pi)$ accounts for permuting the identity into~$\overline\pi$.
\end{proof}

\begin{observation}
  For any ordering~$\sigma$ of~$S$ and any permutation~$\rho$ of the ordered sequence $(i,j,k,l)$,
  the sign of the Pl\"ucker relation only changes with the signature of~$\rho$:
  \[
    \Gamma(\sigma|\rho)
    \ = \
    \varepsilon(\rho)\,
    \Gamma(S|ijkl).
  \]
  % In particular,
  Whether $\Gamma(S|\rho)=0$ or not is independent of~$\epsilon(\rho)$,
  so there is no loss of generality in assuming that
  \[
    1\le i<j<k<l\le n.
  \]
\end{observation}

\section{Pl\"ucker trees}

\subsection{Combining two Pl\"ucker relations}

When trying to prove non-realizability of a simplicial complex~$\Sigma$,
in general we won't be so lucky as to find a positive Pl\"ucker relation for~$\Sigma$ as in Example~\ref{ex:minimal}.
To find a sum of products of known-positive determinants that sums to zero,
we need to combine several relations.

Suppose that two Pl\"ucker relations
$\mathtt{\Gamma_1 = \Gamma(S_1|a_1'a_1''b'b'')}$ and
$\mathtt{\Gamma_2 = \Gamma(S_2|a_2'a_2''c'c'')}$
share the same solid
\[
  \mathtt{
    A
  \ = \
  S_1\cup\{a_1',a_1''\}
  \ = \
  S_2\cup\{a_2',a_2''\},}
\]
possibly with different permutations of its elements.
Written out and partially normalized, they read
\begin{align*}
  \Gamma_1
  \ = \
  \mathtt{\Gamma(S_1|a_1'a_1''b'b'')}
  & \ = \
    \hspace{1.9cm}
    \mathtt{[S_1a_1'a_1'']\,[S_1b'b'']
    {}- [S_1a_1'b']\,[S_1a_1''b''] + [S_1a_1'b'']\,[S_1a_1''b']}
  \\
  & \ = \
    \chi(\mathtt{[A]})\chi(\mathtt{[S_1b'b'']})
    \cdot
    \mathtt{[\overline{A}]}
    \mathtt{[\overline{S_1b'b''}]}
    \mathtt{{}- [S_1a_1'b']\,[S_1a_1''b''] + [S_1a_1'b'']\,[S_1a_1''b']},
    \\[2ex]
  \Gamma_2
  \ = \
  \mathtt{\Gamma(S_2|a_2'a_2''c'c'')}
  & \ = \
    \hspace{1.9cm}
    \mathtt{[S_2a_2'a_2'']\,[S_2c'c''] - [S_2a_2'c']\,[S_2a_2''c''] + [S_2a_2'c'']\,[S_2a_2''c']}
  \\ & \ = \
    \chi(\mathtt{[A]})\chi(\mathtt{[S_2c'c'']})
    \cdot
    \mathtt{[\overline{A}]}
    \mathtt{[\overline{S_2c'c''}]}
       \mathtt{{}- [S_2a_2'c']\,[S_2a_2''c''] + [S_2a_2'c'']\,[S_2a_2''c']}.
\end{align*}
To eliminate the term containing $\mathtt A$ from them, we form the combination
\begin{align}
  \mathtt{
  [\overline{S_2c'c''}]\ \Gamma_1 \ + \
  }s\
  \mathtt{
  [\overline{S_1b'b''}]\ \Gamma_2
  }
  & \ = \
  \mathtt{\phantom{{}+s\ }[\overline{S_2c'c''}]\,\ \big(
    -[S_1a_1'b']\,[S_1a_1''b''] + [S_1a_1'b'']\,[S_1a_1''b']\,
    \big)}
    \label{eq:comb}
  \\
  & \ \phantom{{}={}}\
    {}+s\ \mathtt{\overline{[S_1b'b''}]\,\ \big(
    {}-[S_2a_2'c']\,[S_2a_2''c''] + [S_2a_2'c'']\,[S_2a_2''c']\,\big)},
    \notag
\end{align}
where the \defn{elimination sign} $s$ equals
\[
  s
  \ = \
  -\chi^2(\mathtt{[A]})\,
  \chi(\mathtt{[S_1b'b'']})\,
  \chi(\mathtt{[S_2c'c'']})
  \ = \
  -\chi(\mathtt{[S_1b'b'']})\,
  \chi(\mathtt{[S_2c'c'']}).
\]
In particular, the elimination sign does not depend on $\mathtt{[A]}$.
We represent this elimination as follows:

%in Figure~\ref{fig:elim1}.
%\begin{figure}[htbp]
%  \centering

\ \hfill
\begin{tikzpicture}
    \begin{scope}[
      every node/.style={blue,rectangle,draw,rounded corners=.8ex}
      ]
      \node (G1) at (0,0) { $\mathtt{\Gamma(S_1|a_1'a_1''b'b'')}$ };
      \node (G2) at (4,0) { $\mathtt{\Gamma(S_2|a_2'a_2''c'c'')}$ };
    \end{scope}
    \begin{scope}[
            every node/.style={red, circle, fill=white}
      % every edge./style={very thick}
      ]
      \path[very thick] (G1) edge node{$\TA$} (G2);
    \end{scope}
  \end{tikzpicture}
%  \caption{Eliminating the solid $\mathtt{A}$ between $\mathtt{\Gamma_1}$ and $\mathtt{\Gamma_1}$.}
%  \label{fig:elim1}
%\end{figure}
\hfill\
  
The case of most interest to us is when $\chi(\mathtt{[A]})$ is unknown and $\chi(\mathtt{[S_1b'b'']}) \chi(\mathtt{[S_2c'c'']})= -1$,
because then the elimination sign $s=1$,
and the combination~\eqref{eq:comb} contains one fewer solid of unknown sign.

\subsection{Combining eliminations into trees}

\begin{definition}
  Let $\Sigma$ be a simplicial complex on $n$~vertices of dimension $d-1$.

  \begin{enumerate}
  \item
    The \defn{certificate ring $R[\Sigma]$} is the polynomial ring
    on the $\binom{n}{d+1}$ normal forms $\overline\pi$ of solids
    $\pi\in\binom{[n]}{d+1}$ with coefficients in some base ring~$R$;
    for definiteness, we use $R=\ZZ$.
    In particular, Pl\"ucker relations of~$\Sigma$ lie in~$R[\Sigma]$.
  
  \item
    A \defn{Pl\"ucker tree $T$} for~$\Sigma$ is a tree whose nodes are labeled with Pl\"ucker relations for~$\Sigma$,
    and where each edge is labeled with a solid contained in both relations attached to the incident nodes.
    We will usually confuse nodes and edges with their labels.
    % Two conditions must be fulfilled:
    % \begin{itemize}
    % \item
    %   We require the edge labeling to be injective, so that no solid
    %   appears twice among the edges of~$T$.
    %   
    % \item
      For any two nodes $\mathtt{\Gamma_1}$, $\mathtt{\Gamma_2}$ connected by an edge $\TA$, we require that the monomial containing $\TA$ have
      different sign in $\mathtt{\Gamma_1}$ and $\mathtt{\Gamma_2}$.
%    \end{itemize}

  \item
    \defn{Eliminating an edge} in a Pl\"ucker tree consists of replacing an edge $\TA$
    and its incident nodes $\mathtt{\Gamma_1}$, $\mathtt{\Gamma_2}$
    by a new node containing the combined relation according to~\eqref{eq:comb},
    and connecting it to all neighbors of $\mathtt{\Gamma_1}$, $\mathtt{\Gamma_2}$.
    The labels of the new edges are assigned according to the connected nodes.

  \item
    The \defn{Pl\"ucker certificate $Z(T)$} obtained from~$T$ is the polynomial in~$R[\Sigma]$
    contained in the only node left after all edges of~$T$ have been eliminated.

  \item
    A Pl\"ucker tree is \defn{positive} if $Z(T)$ is a sum of products of normal forms known to be positive
    (see Definition~\ref{def:known}), with positive coefficients.
\end{enumerate}
\end{definition}

Obviously, the existence of a positive Pl\"ucker tree implies the non-realizability of~$\Sigma$.

\begin{lemma}
  \label{lem:eval-order}
  The Pl\"ucker certificate $Z(T)$ % is well-defined,
  % that is, it
  does not depend on the order in which the edges are eliminated.
\end{lemma}

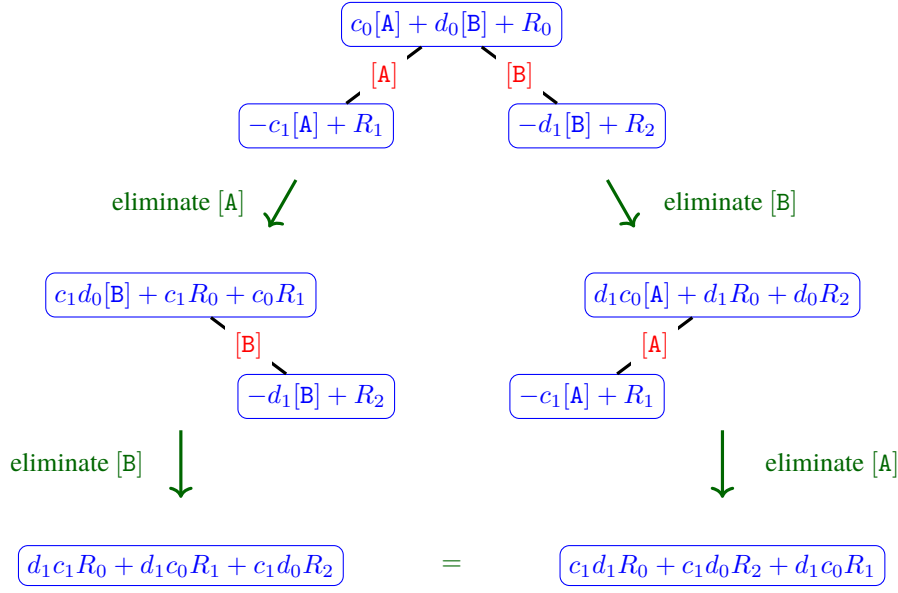
\begin{figure}[htbp]
  \centering
  \begin{tikzpicture}[scale=.9]
    \begin{scope}[
      every node/.style={blue,rectangle,draw,rounded corners=.8ex}
      ]
      \node (G1) at (0,10) { $c_0\TA + d_0\TB + R_0$ };
      \node (G2) at (-2,8.5) { $-c_1\TA + R_1$ };
      \node (G3) at (2,8.5) { $-d_1\TB + R_2$ };

      \node (G4) at (-4,6) { $c_1d_0\TB + c_1 R_0 + c_0R_1$ };
      \node (G5) at (-2,4.5) { $-d_1\TB + R_2$ };      

      \node (G6) at (4,6) { $d_1c_0\TA + d_1 R_0 + d_0R_2$ };
      \node (G7) at (2,4.5) { $-c_1\TA + R_1$ };      

      \node (G8) at (-4,2) { $d_1 c_1 R_0 + d_1 c_0 R_1 +  c_1 d_0 R_2$ };
      \node (G9) at (4,2) { $c_1 d_1 R_0 + c_1 d_0 R_2 + d_1 c_0 R_1 $ };
    \end{scope}
    \begin{scope}[
            every node/.style={red, rectangle, fill=white},
            every path/.style={very thick}
      ]
      \path (G1) edge node{$\TA$} (G2);
      \path (G1) edge node{$\TB$} (G3);
      \path (G4) edge node{$\TB$} (G5);
      \path (G6) edge node{$\TA$} (G7);
    \end{scope}
    \begin{scope}[
      every path/.style={very thick, green!40!black, ->}
      ]
      \path (-2.3,7.7) edge node[left] {eliminate $\TA$\quad\ } (-2.7,7);
      \path (2.3,7.7) edge node[right] {\quad\ eliminate $\TB$} (2.7,7);
      \path (-4,4) edge node[left] {eliminate $\TB$\quad\ } (-4,3);
      \path (4,4) edge node[right] {\quad\ eliminate $\TA$\quad\ } (4,3);
      \node at (0,2) { $=$ };
    \end{scope}
  \end{tikzpicture}  
  \caption{Eliminating edges in any order gives the same result}
  \label{fig:leaves}
\end{figure}

\begin{proof}
  Root $T$ from any node by directing its edges away from that node,
  and consider any node~$N$ of~$T$.
  If $N$~is a leaf, there is nothing to prove.
  If $N$~has exactly one child and this child is a leaf, the corresponding edge may be eliminated uniquely,
  so again there is nothing to prove.

  We may therefore suppose that $N$ has at least two children that are leaves.
  By induction, 
  assume that $N$ is labeled by a polynomial $p\in R[\Sigma]$.
  Further assume that two of its incident leaf nodes are labeled by Pl\"ucker polynomials $\mathtt{\Gamma_1}$ and  $\mathtt{\Gamma_2}$,
  and that the edges $\{p,\mathtt{\Gamma_1}\}$ and $\{p,\mathtt{\Gamma_2}\}$ are labeled by $\mathtt{[A]}$ and $\mathtt{[B]}$, respectively.
  Since the edge labeling is injective, we can moreover assume that
  \begin{itemize}
  \item  $p=c_0\TA + d_0\TB + R_0$ for some
    polynomials $c_0,d_0,R_0\in R[\Sigma]$,
  \item $\mathtt{\Gamma_1} = -c_1\TA + R_1$ for some $c_1,R_1\in R[\Sigma]$, 
  \item $\mathtt{\Gamma_2} = -d_1\TB + R_2$ for some $d_1,R_2\in R[\Sigma]$.
\end{itemize}
  The proof now follows from Figure~\ref{fig:leaves}.
\end{proof}

\begin{remark}
  What we are doing is reminiscent of the approach in Corollary~3.13  in~\cite{bokowski-sturmfels-89},
  except that we only use Pl\"ucker polynomials and not the \emph{van der Waerden identity},
  and we only use monomials and not binomials as coefficients.
\end{remark}

% Solids: $\binom{[n]}{d+1}$, Pl\"ucker relations: $G_0 = \{+,-\}\times\binom{[n]}{d-1}\times\binom{[n-d+1]}{4}$

\subsection{The GP graph: where positive Pl\"ucker trees hide}

Let us define the multigraph on the set of Pl\"ucker relations on which we will search for our Pl\"ucker tree certificates.

 \begin{definition}[Canonical signs; no two adjacent unknown signs] %; $N(\Gamma)$]
   Let $\Gamma=\pm\Gamma(S|ijkl)$ be a Pl\"ucker relation or its negative.
   \begin{enumerate}
   \item The \defn{canonical signs}
     $\sigma_1, \sigma_2, \sigma_3\in\{+,-,?\}$ of $\Gamma$ are
   \begin{align*}
     \bblue{\sigma_1}
     \ := \
     \pm\chi_{[Sij]}\chi_{[Skl]},
     % \\
     \qquad
     \bblue{\sigma_2}
     \ := \
     \mp \chi_{[Sik]}\chi_{[Sjl]},
     \qquad
     \bblue{\sigma_3}
     \ := \
       \pm
     \chi_{[Sil]}\chi_{[Sjk]},
   \end{align*}
   so that  
   \begin{align*}
     \pm\Gamma(S|ijkl)
     &\ = \
       \pm[Sij][Skl] 
       \mp [Sik][Sjl]
       \pm [Sil][Sjk]
     \\
     &\ = \
       \pm\chi_{[Sij]}\chi_{[Skl]}\overline{[Sij]}\ \overline{[Skl]}
       \mp \chi_{[Sik]}\chi_{[Sjl]}\overline{[Sik]}\ \overline{[Sjl]}
       \pm \chi_{[Sil]}\chi_{[Sjk]}\overline{[Sil]}\ \overline{[Sjk]}
     \\
     &\ =: \
       \bblue{\sigma_1}\overline{[Sij]}\ \overline{[Skl]}
       +\bblue{\sigma_2} \overline{[Sik]}\ \overline{[Sjl]}
       +\bblue{\sigma_3}\overline{[Sil]}\ \overline{[Sjk]}
     \\
     &\ =: \
       \sigma_1 \blue{s_{11}s_{12}}
       +\sigma_2 \blue{s_{21}s_{22}}
       +\sigma_3 \blue{s_{31}s_{32}}
       .
   \end{align*}
   Here we multiply unknown signs via the rules $(+)(?) =\ ?$ and $(-)(?) =\ ?$.
   This is called a \emph{fuzzy ring} in~\cite[Chapter~3]{bokowski-sturmfels-89}.

 \item
   The relation $\Gamma$ has \defn{no two adjacent unknown solids} if for each $i=1,2,3$,
   at least one of the signs of the solids corresponding to~$\sigma_i$ is known.

 % \item
 %   The set of \defn{negative solids} of~$\Gamma$ is $\bblue{N(\Gamma)}=\big\{\{s_{i1},s_{i2}\} : \sigma_i=-\big\}$.
 \end{enumerate}
\end{definition}

Why do we not want a relation $\Gamma$ to have two adjacent unknown solids? Look at~\eqref{eq:comb} and~Figure~\ref{fig:leaves}:
we use one of the solids as a coefficient for elimination, so it will most likely survive into the final certificate;
therefore, we had better know its orientation.

\begin{definition}[GP graph]
  Let $\Sigma$ be a simplicial complex,
  and let $\overline S_\Sigma$ be the set of canonical solids in~$\Sigma$.
  The \defn{GP graph $\mGP(\Sigma)$}${}=(V_\Sigma,E_\Sigma)$ of~$\Sigma$ is an undirected graph %$(V_\Sigma,E_\Sigma)$ % with set of undirected colored edges
  % $E_\Sigma \subseteq \overline S_\Sigma\times \binom{V_\Sigma}{2}$
  on the node set
  \[
    \bblue{V_\Sigma}
    \ = \
    \Big\{
    \pm\Gamma(S|ijkl) :
    % \sigma_1,\sigma_2,\sigma_3\in\{+,?\}
    \Gamma(S|ijk\ell)
    \text{ has no two adjacent unknown solids }
    \Big\}.
  \]
  The edges of $\mGP$ are colored with the set $\mathfrak{S}(\Sigma)$ of normal forms of solids of~$\Sigma$.
  There can be multiple edges between two nodes, but they will have different colors.
  More precisely, two nodes $\Gamma,\Gamma'\in V$  are joined with an edge $(\overline\pi,\{\Gamma,\Gamma'\})$
  colored~$\overline{\pi}\in\mathfrak{S}(\Sigma)$ in~$\bblue{E_\Sigma}\subseteq \mathfrak{S}(\Sigma)\times \binom{V_\Sigma}{2}$ iff they share a solid~$\pi$ such that
  \[
    \sigma_i
    \ = \
    -\sigma_i',
  \]
  where $\sigma_i,\sigma_i'$ are the canonical signs of the terms containing $\pi$ in~$\Gamma,\Gamma'$.
\end{definition}

\begin{example}[Multiple edges in $\mGP(\Sigma)$]
  Suppose that 
  \begin{align*}
    \Gamma(S|xjkl)
    & \ = \
      [Sxj][Skl] - [Sxk][Sjl] + [Sxl][Sjk]
      \ = \
      \phantom{+}\overline{[Sxj]}\,\overline{[Skl]} 
      + \overline{[Sxk]}\,\overline{[Sjl]}
      + \overline{[Sxl]}\,\overline{[Sjk]}
      ,
    \\
    \Gamma(S|yjkl)
    & \ = \
      [Syj][Skl] - [Syk][Sjl] + [Syl][Sjk]
      \ = \
      - \overline{[Syj]}\,\overline{[Skl]} 
      - \overline{[Syk]}\,\overline{[Sjl]}
      - \overline{[Syl]}\,\overline{[Sjk]}
  \end{align*}
  are nodes of~$\mGP(\Sigma)$.
  Then they
  are connected by three edges with colors $\overline{[Skl]}$, $\overline{[Sjl]}$ and $\overline{[Sjk]}$.
\end{example}

A Pl\"ucker tree $T$ is an induced subgraph of the GP~graph such that
\begin{itemize}
% \item
%   the edge coloring is injective, i.e., no edge color appears twice;

\item
  for each unknown solid $\bred{\overline\pi^?}$ (i.e., $\chi_\pi=\ ?$) in a node~$\Gamma$,
  there is an edge colored~$\overline\pi$ incident to~$\Gamma$;

\item
  $T$ is a tree, i.e, $|V(T)|=|E(T)|+1$.
\end{itemize}

\subsection{Using integer programming}
\label{sec:ip}

We search for such a tree by solving the integer program on the integer indicator variables
$\{x_\Gamma : \Gamma\in V_\Sigma\}$ and
$\{x_e : e=(\overline\pi,\{\Gamma,\Gamma'\})\in E_\Sigma\}$
defined as follows:
\begin{align}
  \min \sum_{\Gamma\in V_\Sigma} x_\Gamma
  \notag
  &
  \\
  \text{s.t.} \hspace{1cm}
    2 \sum_{\overline\pi:\ e=(\overline\pi,\{\Gamma,\Gamma'\})\in E_\Sigma} x_e
  &\ \le \
    x_\Gamma + x_{\Gamma'}
    \qquad\text{for each }
    \{\Gamma,\Gamma'\}\in\binom{V_\Sigma}{2}
    \label{ineq:ggc}
  \\
  \sum_{\overline\pi:\ e=(\overline\pi,\{\Gamma,\Gamma'\})\in E_\Sigma}
  x_e
  &\ \le \
    1
    \qquad\qquad\quad\
    \text{for each }\{\Gamma,\Gamma'\}\in\binom{V_\Sigma}{2},
    \label{ineq:connect}
  \\
    1 + 
    \sum_{e\in E_\Sigma} x_e
  &\ = \
    \sum_{\Gamma\in V_\Sigma} x_\Gamma
    \label{eq:tree}
  \\
  \sum_{\Gamma' :\ e=(\bred{\overline\pi^?},\{\Gamma,\Gamma'\})\in E_\Sigma}
  x_e
  &\ = \
    x_\Gamma
    \qquad
    \text{ for all }\Gamma\in V_\Sigma,
    \text{for all unknown } \bred{\overline\pi^?}\in\Gamma
    \label{ineq:unknown}
  % \\
  % \sum_{j=1,2}\
  % \sum_{\Gamma' :\ e=(\overline\pi_j,\{\Gamma,\Gamma'\})\in E_\Sigma}
  % x_e
  % &\ = \
  %   x_\Gamma
  %   \qquad
  %   \text{ for all }\Gamma\in V_\Sigma,
  %   \text{for all  } \{\overline\pi_1,\overline\pi_2\}\in N(\Gamma),
  %   \label{ineq:negative}
  % \\
  % x_\Gamma,\ x_e
  % & \ \in \
  %   \{0,1\}.
  %   \notag
 \end{align}
 These inequalities for
   $x_\Gamma,\ x_e
   \in 
    \{0,1\}$
have the following interpretation:
 \begin{itemize}
 \item
   \eqref{ineq:ggc} ensures that both endpoints of an edge present in the solution are present;

 \item
   \eqref{ineq:connect} ensures that at most one edge is selected between two selected nodes;

 \item
   \eqref{eq:tree} forces the solution to be a tree with at least one node;

 \item
   \eqref{ineq:unknown} ensures that if a node $\Gamma$~with an unknown sign $\bred{\overline\pi^?}$ is present in the solution,
   then there is exactly one edge of that color incident to~$\Gamma$;

 % \item
 %   \eqref{ineq:negative} ensures that if a node $\Gamma$~with an negative sign is present in the solution,
 %   then there is exactly one edge of an associated solid $\overline\pi_1,\overline\pi_2$ incident to~$\Gamma$.
   
 \end{itemize}

\section{Results}

We have implemented a search for positive Pl\"ucker trees in the software framework
\texttt{polymake}~\cite{polymake:2000}. % where it is available from version 4.1.1 onwards.
Obviously, there are many details of this implementation which we will not discuss in this paper,
such as taking into account possible symmetries of the simplicial complex,
using appropriate data structures, etc.
%\todo{When done, make the source code available at \url{https://polymake.org}.}

We remark that in all cases, the corresponding integer program was solved to optimality,
proving the exhibited positive Pl\"ucker trees to be smallest possible.

\subsection{Zheng's balanced $3$-sphere}

  In \cite{zheng20}, Hailun Zheng constructs a combinatorial $3$-sphere~$Z$ on $16$~vertices with $f$-vector $(16,96,160,80)$ and an action of the dihedral group~$D_4$ of order~$8$.
  The $1$-skeleton of~$Z$ is the complete multipartite graph~$K_{4,4,4,4}$, which implies that it is \emph{balanced 2-neighborly}.

  To explain why this sphere is important, let's fix definitions.
  A $(d-1)$-dimensional simplicial complex is \defn{balanced} if its $1$-skeleton is \defn{$d$-colorable} in the graph-theoretic sense,
  i.e., its vertices can be colored with $d$~colors in such a way that the endpoints of all edges receive different colors.
  Moreover, a $(d-1)$-dimensional balanced simplicial complex~$\Sigma$ is \defn{balanced $k$-neighborly}
  if each $k$-subset of the vertex set that contains at most one vertex of each color class is actually a face of~$\Sigma$.

  Now we can say why Zheng's example is important~---~in fact, it is important in at least \emph{two} ways.
  
  First, there has been a lot of work on analogies between combinatorial data in the balanced and the non-balanced settings~\cite{kub18,kub18a,venturello2019}.
  For example, one would like to have a balanced analogue of the celebrated Upper Bound Theorem by McMullen and Stanley.
  For this, in particular one would like balanced analogues of the extremal examples to even \emph{exist}, i.e.,
  one would like to construct infinite families of balanced $k$-neighborly simplicial spheres.
  What Zheng shows in~\cite{zheng20}, however, is that
  (i) there is no balanced 2-neighborly homology 3-sphere on 12 vertices;
  (ii) there is no balanced 2-neighborly homology 4-sphere on 15 vertices;
  (iii) but taking suspensions over her example~$Z$ yields a balanced 2-neighborly homology $(3 + m)$-sphere on $16 + 2m$ vertices for every $m\ge0$.

  The second reason why her example is important lies in the fact that in \cite{pz12},
  the present author, Vincent Pilaud and Francisco Santos study which graphs are realizable as the $1$-skeleton of polytopes.
  The case of multipartite graphs was not treated there,
  and to date the only polytope whose graph is known to be the multipartite graph $K_{4,4,4,4}$ is the $4$-dimensional cross polytope.

  We can now show for the first time that $Z$~is not realizable as the boundary complex of a convex polytope,
  and therefore that $Z$~does not yield a new polytope whose graph is~$K_{4,4,4,4}$.

  \begin{theorem}
    Zheng's balanced sphere $Z$ is not polytopal.
  \end{theorem}

  \begin{proof}
  We enumerate and orient the facets of~$Z$:
  \begin{quote}\ttfamily\small
   -[048c] +[048e] +[049c] -[049d] +[04ad] -[04ae] +[059d] -[059f] -[05ad] +[05ae] -[05be] +[05bf] +[068c] -[068e] -[069c] +[069e] -[079e] +[079f] +[07be] -[07bf] +[148c] -[148e] +[14ae] -[14af] -[14bc] +[14bf] -[158c] +[158d] -[159d] +[159f] +[15bc] -[15bf] +[168e] -[168f] -[16ae] +[16af] -[178d] +[178f] +[179d] -[179f] -[24ad] +[24af] +[24bd] -[24bf] +[258c] -[258d] -[25ac] +[25ad] -[268c] +[268d] +[269c] -[269e] +[26ae] -[26af] -[26bd] +[26bf] -[279c] +[279e] +[27ac] -[27ae] -[349c] +[349d] +[34bc] -[34bd] +[35ac] -[35ae] -[35bc] +[35be] -[368d] +[368f] +[36bd] -[36bf] +[378d] -[378f] +[379c] -[379d] -[37ac] +[37ae] -[37be] +[37bf] 
  \end{quote}

\noindent The non-realizability follows from the existence of the following Pl\"ucker tree:
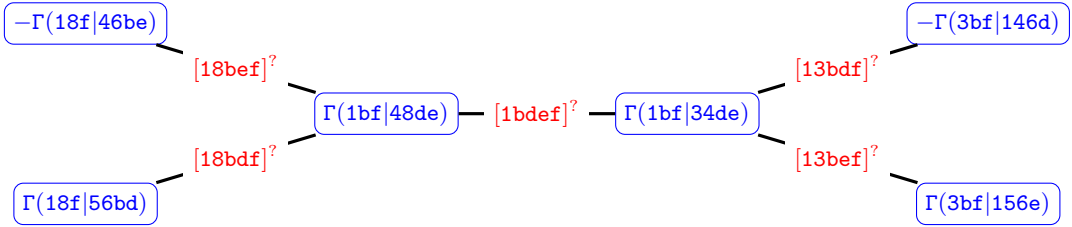
\begin{center}
  \begin{tikzpicture}[scale=.8]
% 18f46ben  0
% 18f56bdp  1
% 1bf34dep  2
% 1bf48dep  3
% 3bf146dn  4
% 3bf156ep  5
% [13bdf]_1bf34dep_3bf146dn  2 4
% [13bef]_1bf34dep_3bf156ep  2 5
% [18bdf]_18f56bdp_1bf48dep  1 3
% [18bef]_18f46ben_1bf48dep  0 3 
% [1bdef]_1bf34dep_1bf48dep  2 3    
    \begin{scope}[
      every node/.style={blue,rectangle,draw,rounded corners=.8ex}
      ]
      \node (G0) at (-7.5,1.5) {\small$-\mathtt{\Gamma(18f|46be)}$ };
      \node (G1) at (-7.5,-1.5) {\small$\mathtt{\Gamma(18f|56bd)}$ };
      \node (G2) at (2.5,0) {\small$\mathtt{\Gamma(1bf|34de)}$ };
      \node (G3) at (-2.5,0) {\small$\mathtt{\Gamma(1bf|48de)}$ };
      \node (G4) at (7.5,1.5) {\small$-\mathtt{\Gamma(3bf|146d)}$ };
      \node (G5) at (7.5,-1.5) {\small$\mathtt{\Gamma(3bf|156e)}$ };

    \end{scope}
    \begin{scope}[
            every node/.style={red, rectangle, fill=white},
            every path/.style={very thick}
      ]
      \path (G2) edge node{\small$\mathtt{[13bdf]^?}$} (G4);
      \path (G2) edge node{\small$\mathtt{[13bef]^?}$} (G5);
      \path (G1) edge node{\small$\mathtt{[18bdf]^?}$} (G3);
      \path (G0) edge node{\small$\mathtt{[18bef]^?}$} (G3);
      \path (G2) edge node{\small$\mathtt{[1bdef]^?}$} (G3);
    \end{scope}
  \end{tikzpicture}    
\end{center}

which gives rise to the following certificate, which is short enough to write down in full:
\small
  \begin{align*}
     \quad
      \mathtt{
      [36fb5] \Bigg([36fb4] \bigg([14bf3] \Big([16f85] \big([14bfd] (-\bblue{\Gamma(18f|46be)})
         } &
    \\ 
         \mathtt{
    + [16f84] \,\bblue{\Gamma(1bf|48de)}\big)
         } &
    \\ 
         \mathtt{
    + [16f84]\,[14bfe]\, \bblue{\Gamma(18f|56bd)}\Big)
      } &
    \\ 
         \mathtt{
     + [16f84]\,[14bf8]\,[16f85] \, \bblue{\Gamma(1bf|34de)}\bigg)
      } &
    \\ 
    \mathtt{
    + [16f84]\,[14bf8]\,[14bfe]\,[16f85] (-\bblue{\Gamma(3bf|146d)})\Bigg)
    } &
    \\ 
    \mathtt{
    \ {}+ [16f84]\,[14bf8]\,[14bfd]\,[16f85]\,[36fb4] \,\bblue{\Gamma(3bf|156e)}\quad
    } \rlap{\ .} &
  \end{align*}
  \normalsize
  Substituting the appropriately normalized Pl\"ucker polynomials yields
  \small
  \begin{align*}
  \mathtt{
    [36fb5] \Bigg([36fb4] \bigg([14bf3] \Big([16f85] \big([14bfd] (-[16f84]\bred{[18bef]^?}+[14bf8][16f8e]+[14e8f][16f8b])
    } &
    \\ 
         \mathtt{
         {}+ [16f84] ([14bf8]\bred{[1bdef]^?}+[14bfd]\bred{[18bef]^?}-[14bfe]\bred{[18bdf]^?})\big)
         } &
    \\ 
         \mathtt{
    {} + [16f84]\,[14bfe] \big([16f85]\bred{[18bdf]^?}+[15fb8][16f8d]+[158df][16f8b]\big)\Big)
         } &
    \\ 
         \mathtt{
    {} + [16f84]\,[14bf8]\,[16f85] \big(-[14bf3]\bred{[1bdef]^?}+\bred{[13bdf]^?}[14bfe]-\bred{[13bef]^?}[14bfd]\big)\bigg)
         } &
    \\ 
         \mathtt{
    {} + [16f84]\,[14bf8]\,[14bfe]\,[16f85] \big([14bf3][36fbd]+[36fb1][34dbf]-\bred{[13bdf]^?}[36fb4]\big)\Bigg)
         } &
    \\ 
         \mathtt{
    {} + [16f84]\,[14bf8]\,[14bfd]\,[16f85]\,[36fb4] \big([15fb3][36fbe]+[36fb1][35bef]+\bred{[13bef]^?}[36fb5]\big)
    }\rlap{\ ,}
\end{align*}
\normalsize
and this in turn simplifies to the final form of the certificate,
\small
\begin{align*}
  \mathtt{
    [14bf3]\,[16f84]\,[14bf8]\,[14bfe]\,[16f85]\,[36fb5]\,[36fbd] + [14bf3]\,[16f84]\,[14bfe]\,[15fb8]\,[16f8d]\,[36fb4]\,[36fb5]
         } &
    \\ 
         \mathtt{
    {} + [14bf3]\,[16f84]\,[14bfe]\,[158df]\,[16f8b]\,[36fb4]\,[36fb5] + [14bf3]\,[14bf8]\,[14bfd]\,[16f85]\,[16f8e]\,[36fb4]\,[36fb5]
         } &
    \\ 
         \mathtt{
    {} + [14bf3]\,[14e8f]\,[14bfd]\,[16f85]\,[16f8b]\,[36fb4]\,[36fb5] + [15fb3]\,[16f84]\,[14bf8]\,[14bfd]\,[16f85]\,[36fb4]\,[36fbe]
         } &
    \\ 
         \mathtt{
    {} + [36fb1]\,[16f84]\,[14bf8]\,[14bfd]\,[16f85]\,[36fb4]\,[35bef] + [36fb1]\,[16f84]\,[14bf8]\,[14bfe]\,[16f85]\,[34dbf]\,[36fb5]
    \rlap{\ .}
    } &
  \end{align*}
\normalsize
\end{proof}

\subsection{Topological Prismatoids}

In~\cite{prismatoids2019}, Francisco Criado and Francisco Santos introduced \emph{topological prismatoids},
a combinatorial abstraction of the geometric prismatoids used by Santos~\cite{santos2012} to construct counterexamples to the Hirsch conjecture.
Criado and Santos construct four combinatorially distinct non-$d$-step topological $4$-dimensional prismatoids on $14$~vertices,
referred to as \#1039, \#1963, \#2669 and \#3513,
which imply the existence of $8$-dimensional spheres on $18$~vertices whose combinatorial diameter exceeds the Hirsch bound.
In~\cite{prismatoids2019}, the question of polytopality of these combinatorial prismatoids was left open.
Using our technique, we can prove that all four of them are in fact not polytopal.

\begin{remark}
  During the elaboration of this manuscript, Gouveia, Macchia and Wiebe (see version \texttt{v3} of~\cite{gouveia-macchia-wiebe20})
  were also able to prove the non-realizability of the prismatoid \#1039;
  they will make non-realizability proofs of the other ones available in a future paper~\cite{wiebe-personal}.
\end{remark}

\begin{remark}
A technical detail worth mentioning is that each of these combinatorial prismatoids $\Pi$ has two non-simplicial facets.
To construct an orientation of the boundary of~$\Pi$, we first consider the simplicial sphere~$\widetilde\Pi$
obtained by replacing each non-simplicial facet by a cone over its (simplicial) boundary.
Next, we calculate a representative cycle for the top simplicial homology group of~$\widetilde\Pi$,
and restrict that cycle to the simplicial facets of~$\Pi$.
\end{remark}

\begin{theorem}
  \label{thm:1039}
   The prismatoid \#1039 is not realizable.
 \end{theorem}

 \begin{proof}
   One of the two possible orientations of the simplicial facets of \#1039 is as follows.
   The nine facets relevant to the certificate are listed first:
   \begin{quote}\ttfamily\small
     \noindent
     \bblue{-[0145f] +[014ad] -[014ae] -[014bc] +[014bf] -[014cd] -[04bcd] -[05bde] +[05cde]}

     \medskip
     \noindent
     -[0123d] +[0126d] -[0134e] -[013ad] +[013ae] -[0156g] -[015bf] +[015bg] +[016cd] -[016cg] +[01bcg] +[0234e] +[023cd] -[023ce] +[0245f] +[024ae] -[024af] +[0256g] -[025ae] +[025af] +[025be] -[025bg] -[026be] +[026bg] -[026cd] +[026ce] +[03acd] -[03ace] -[04abd] +[04abf] +[05abd] -[05abf] -[05acd] +[05ace] +[06bce] -[06bcg] -[0bcde] -[1234e] +[123ae] -[123af] +[123bf] -[123bg] -[123cd] +[123cg] -[1245f] -[124ae] +[124af] -[1256g] -[125bf] +[125bg] +[126cd] -[126cg] +[13abf] -[13abg] -[13acd] +[13acg] -[14abf] +[14abg] +[14acd] -[14acg] +[14bcg] +[23ace] -[23acf] +[23bcf] -[23bcg] -[25abe] +[25abf] +[26abe] -[26abf] -[26ace] +[26acf] -[26bcf] +[26bcg] -[3abfg] +[3acfg] -[3bcfg] -[4abcd] +[4abcg] -[5abde] +[5acde] -[6abef] +[6acef] -[6bcef] 
   \end{quote}

The non-realizability follows from the existence of the following positive Pl\"ucker tree:

% node 1519: <0 1 4 5 10 13 14 15> 0 = 0145adefp
% node 2178: <0 1 4 13 5 11 12 15> 0 = 014d5bcfp
% node 12498: <0 4 5 13 1 11 12 14> 0 = 045d1bcep
% edge 307 (12498,1519): [0145de]_045d1bcep_0145adefp
% edge 299 (12498,2178): [0145bd]_045d1bcep_014d5bcfp

% \Gamma_0 = \bblue{\Gamma(0145|adef)}
% \Gamma_1 = \bblue{\Gamma(014d|5bcf)}
% \Gamma_2 = \bblue{\Gamma(045d|1bce)}
% tree connectivity:
% {2}
% {2}
% {0 1}

  {\hfill%\begin{center}
  \begin{tikzpicture}
    \begin{scope}[
      every node/.style={blue,rectangle,draw,rounded corners=.8ex}
      ]
      \node (G0) at (0,0) {\small$\mathtt{\Gamma(0145|adef)}$ };
      \node (G1) at (10,0) {\small$\mathtt{\Gamma(014d|5bcf)}$ };
      \node (G2) at (5,0) {\small$\mathtt{\Gamma(045d|1bce)}$ };

    \end{scope}
    \begin{scope}[
            every node/.style={red, circle, fill=white},
            every path/.style={very thick}
      ]
      \path (G0) edge node{\small$\mathtt{[0145de]^?}$} (G2);
      \path (G1) edge node{\small$\mathtt{[0145bd]^?}$} (G2);
    \end{scope}
  \end{tikzpicture}    
\hfill%\end{center}
}

which gives rise to the certificate
\small
\begin{align*}
  \mathtt{
  [014ad5]\,[014f5e]\,[014dcf]\,[04bdc5] + [014ea5]\,[014f5d]\,[014dcf]\,[04bdc5]
  }
  & \\
  \mathtt{
  {}+ [014f5a]\,[014dc5]\,[014bfd]\,[05cde4] + [014f5a]\,[014dc5]\,[014dcf]\,[05bed4]
  }
  & \\
  \mathtt{
  {}+ [014f5a]\,[014f5d]\,[04bdc1]\,[05cde4]
  \rlap{\, .}
  }
\end{align*}
\normalsize
\end{proof}

\begin{theorem}
  \label{thm:1963}
  The prismatoid \#1963 is not realizable.
\end{theorem}

\begin{proof}
  The oriented simplicial facets are as follows, relevant ones first:
  \begin{quote}\ttfamily\small
    \noindent
    \bblue{-[015cf] +[0245f] +[0256g] +[025ce] -[025cg] -[06bcf] -[125cf] -[25bce] }

    \medskip
    \noindent
    -[0123d] +[0126d] -[0134e] +[013ae] -[013af] -[013bd] +[013bg] +[013cf] -[013cg] -[0145f] -[014ae] +[014af] -[0156g] +[015cg] +[016bd] -[016bg] +[0234e] +[023bd] -[023be] +[024ae] -[024af] -[025ae] +[025af] -[026bd] +[026be] -[026ce] +[026cg] -[03abe] +[03abf] +[03bcf] -[03bcg] +[05ace] -[05acf] +[06abe] -[06abf] -[06ace] +[06acf] +[06bcg] -[1234e] -[123bd] +[123be] -[1245f] -[124ac] +[124af] +[124bd] -[124be] -[124cd] -[1256g] +[125cg] +[126cd] -[126cg] +[12acf] +[13abe] -[13abg] -[13acf] +[13acg] -[14abe] +[14abg] -[14acg] -[14bcd] +[14bcg] +[16bcd] -[16bcg] -[24abd] +[24abe] +[24acd] +[25abd] -[25abe] -[25acd] +[25acf] +[25bcd] -[26bcd] +[26bce] -[3abfg] +[3acfg] -[3bcfg] -[4abcd] +[4abcg] -[5abde] +[5acde] -[5bcde] -[6abef] +[6acef] -[6bcef] 
  \end{quote}

  The non-realizability follows from the existence of the following positive Pl\"ucker tree:
  % node 8109: <0 2 5 15 4 6 12 16> 1 = 025f46cgn
  % node 15599: <0 5 12 15 1 2 6 11> 1 = 05cf126bn
  % node 38774: <2 5 12 15 0 1 11 14> 0 = 25cf01bep
  % edge 712 (15599,8109): [0256cf]_05cf126bn_025f46cgn
  % edge 737 (38774,15599): [025bcf]_25cf01bep_05cf126bn

  % \Gamma_0 = \bblue{(-\Gamma(025f|46cg))}
  % \Gamma_1 = \bblue{(-\Gamma(05cf|126b))}
  % \Gamma_2 = \bblue{\Gamma(25cf|01be)}
  % tree connectivity:
  % {1}
  % {0 2}
  % {1}

  \hfill
  \begin{tikzpicture}
    \begin{scope}[
      every node/.style={blue,rectangle,draw,rounded corners=.8ex}
      ]
      \node (G1) at (0,0) { $\mathtt{-\Gamma(025f|46cg)}$ };
      \node (G2) at (5,0) { $\mathtt{-\Gamma(05cf|126b)}$ };
      \node (G3) at (10,0) { $\mathtt{\Gamma(25cf|01be)}$ };

    \end{scope}
    \begin{scope}[
      every node/.style={red, circle, fill=white},
      every path/.style={very thick}
      ]
      \path (G1) edge node{$\mathtt{[0256cf]^?}$} (G2);
      \path (G2) edge node{$\mathtt{[025bcf]^?}$} (G3);
    \end{scope}
  \end{tikzpicture}    
  \hfill
  \

  with final certificate
  \small
  \begin{align*}
    \mathtt{
    [125fc0]\,[015fc6]\,[0245fg]\,[25becf] + [125fc0]\,[0245fg]\,[06bfc5]\,[125fce]
    }
    & \\
    \mathtt{
    {}+ [015fc6]\,[0245fg]\,[025cef]\,[125fcb] + [015fcb]\,[0245f6]\,[025gcf]\,[125fce]
    }
    &\\
    \mathtt{
    {}+ [015fcb]\,[0245fc]\,[0256gf]\,[125fce]
    \rlap{\, .}
    }
  \end{align*}
  \normalsize
\end{proof}

\begin{theorem}
  \label{thm:2669}
  The prismatoid \#2669 is not realizable.
\end{theorem}

\begin{proof}
  The oriented simplicial facets of \#2669 are (relevant ones first)
  \begin{quote}\ttfamily\small
    \noindent
    \bblue{+[0234a] +[023ad] +[123ae] -[123cd] -[13acd] -[14abd] +[14acd] -[234ae] +[23abf] +[23acd] -[24abf] -[4abcd]}

    \medskip
    \noindent
    -[0123d] +[0126d] -[0134e] -[013ad] +[013ae] -[0145f] -[014be] +[014bf] -[0156g] +[015ad] -[015ae] +[015be] -[015bf] -[015cd] +[015cg] +[016cd] -[016cg] +[0245f] -[024af] +[0256g] +[025bf] -[025bg] -[026ad] +[026af] -[026bf] +[026bg] +[034ae] -[04abe] +[04abf] -[05acd] +[05ace] -[05bce] +[05bcg] +[06abe] -[06abf] +[06acd] -[06ace] +[06bce] -[06bcg] -[1234e] -[123ag] +[123cg] -[1245f] -[124ae] +[124ag] +[124bf] -[124bg] -[1256g] -[125bf] +[125bg] +[126cd] -[126cg] +[13acg] +[14abe] -[14acg] -[14bcd] +[14bcg] +[15abd] -[15abe] +[15bcd] -[15bcg] -[23abg] -[23acf] +[23bcf] -[23bcg] +[24abg] -[26acd] +[26acf] -[26bcf] +[26bcg] -[3abfg] +[3acfg] -[3bcfg] +[4abcg] -[5abde] +[5acde] -[5bcde] -[6abef] +[6acef] -[6bcef] 
  \end{quote}

  The non-realizability follows from the existence of the following positive Pl\"ucker tree:
  % node 33198: <2 3 10 11 0 4 13 15> 1 = 23ab04dfn
  % node 33610: <2 3 10 13 1 4 12 14> 1 = 23ad14cen
  % node 36114: <2 4 10 13 1 3 11 12> 0 = 24ad13bcp
  % edge 2129 (36114,33198): [234abd]_24ad13bcp_23ab04dfn
  % edge 1304 (36114,33610): [1234ad]_24ad13bcp_23ad14cen

  % \Gamma_0 = \bblue{(-\Gamma(23ab|04df))}
  % \Gamma_1 = \bblue{(-\Gamma(23ad|14ce))}
  % \Gamma_2 = \bblue{\Gamma(24ad|13bc)}
  % tree connectivity:
  % {2}
  % {2}
  % {0 1}

  \hfill
  \begin{tikzpicture}
    \begin{scope}[
      every node/.style={blue,rectangle,draw,rounded corners=.8ex}
      ]
      \node (G0) at (0,0) { $\mathtt{-\Gamma(23ab|04df)}$ };
      \node (G1) at (10,0) { $\mathtt{-\Gamma(23ad|04df)}$ };
      \node (G2) at (5,0) { $\mathtt{\Gamma(24ad|13bc)}$ };

    \end{scope}
    \begin{scope}[
      every node/.style={red, circle, fill=white},
      every path/.style={very thick}
      ]
      \path (G0) edge node{$\mathtt{[234abd]^?}$} (G2);
      \path (G1) edge node{$\mathtt{[1234ad]^?}$} (G2);
    \end{scope}
  \end{tikzpicture}    
\hfill
\

with final certificate
\small
\begin{align*}
  \mathtt{
  [0234ab]\,[14acd2]\,[23abfd]\,[23acde] + [023adb]\,[14acd2]\,[24afb3]\,[23acde]
    }
  & \\
  \mathtt{
  {} + [23abf0]\,[23acd1]\,[234ead]\,[4abdc2] + [23abf0]\,[123aed]\,[23acd4]\,[4abdc2]
    }
  & \\
  \mathtt{
  {}+ [23abf0]\,[14adb2]\,[23acd4]\,[23acde]
  \rlap{\, .}
  }
\end{align*}

 \end{proof}

 \begin{theorem}
   \label{thm:3513}
    The prismatoid \#3513 is not realizable.
  \end{theorem}

  \begin{proof}
    The oriented simplicial facets are (relevant ones first)
    \begin{quote}\ttfamily\small
      \noindent
      \bblue{+[014af] -[014ag] +[014bg] +[04abg] +[124af] +[125bg] -[13abg] -[14acg] +[14bcg] +[15abg] -[3abfg] +[4abcg]}

      \medskip
      \noindent
      -[0123d] +[0126d] -[0134e] -[013cd] +[013ce] -[0145f] -[014bc] -[014ce] -[0156g] -[015af] +[015ag] +[016bd] -[016bg] +[01bcd] +[0234e] +[023cd] -[023ce] +[0245f] +[024ae] -[024af] +[0256g] -[025ae] +[025af] +[025bd] -[025bg] -[025cd] +[025ce] -[026bd] +[026bg] -[04abd] +[04acd] -[04ace] -[04bcd] +[05abd] -[05abg] -[05acd] +[05ace] -[1234e] +[123ae] -[123af] -[123cd] +[123cf] -[1245f] -[124ae] -[1256g] -[125bf] +[126bf] -[126bg] +[126cd] -[126cf] +[13abf] -[13ace] +[13acg] +[13bcf] -[13bcg] +[14ace] -[15abf] +[16bcd] -[16bcf] +[23ace] -[23acf] -[25abe] +[25abf] +[25bcd] -[25bce] +[26abe] -[26abf] -[26ace] +[26acf] -[26bcd] +[26bce] +[3acfg] -[3bcfg] -[4abcd] -[5abde] +[5acde] -[5bcde] -[6abef] +[6acef] -[6bcef] 
    \end{quote}

The non-realizability follows from the existence of the following positive Pl\"ucker tree:

% node 27490: <1 4 10 16 0 2 11 15> 0 = 14ag02bfp
% node 27904: <1 4 11 16 0 2 5 10> 0 = 14bg025ap
% node 46246: <4 10 11 16 1 3 12 15> 1 = 4abg13cfn
% edge 1460 (27904,27490): [124abg]_14bg025ap_14ag02bfp
% edge 1927 (46246,27490): [14abfg]_4abg13cfn_14ag02bfp

% \Gamma_0 = \bblue{\Gamma(14ag|02bf)}
% \Gamma_1 = \bblue{\Gamma(14bg|025a)}
% \Gamma_2 = \bblue{(-\Gamma(4abg|13cf))}
% tree connectivity:
% {1 2}
% {0}
% {0}
\hfill
  \begin{tikzpicture}
    \begin{scope}[
      every node/.style={blue,rectangle,draw,rounded corners=.8ex}
      ]
      \node (G0) at (5,0) { $\mathtt{\Gamma(14ag|02cf)}$ };
      \node (G1) at (0,0) { $\mathtt{\Gamma(14bg|025a)}$ };
      \node (G2) at (10,0) { $-\mathtt{\Gamma(4abg|13cf)}$ };

    \end{scope}
    \begin{scope}[
            every node/.style={red, circle, fill=white},
            every path/.style={very thick}
      ]
      \path (G1) edge node{$\mathtt{[124abg]^?}$} (G0);
      \path (G2) edge node{$\mathtt{[14abfg]^?}$} (G0);
    \end{scope}
  \end{tikzpicture}    
\hfill
\

with final certificate
\small
\begin{align*}
  \mathtt{
   [014ga2]\,[014bg5]\,[13agb4]\,[4abcgf] + [014ga2]\,[014bg5]\,[4abcg1]\,[3abgf4]     }
  & \\
  \mathtt{
  {} + [014bg2]\,[014gaf]\,[15abg4]\,[4abcg3] + [014bg5]\,[04abg1]\,[124afg]\,[4abcg3]     }
  & \\
  \mathtt{
  {} + [04abg1]\,[014gaf]\,[125bg4]\,[4abcg3]
  \rlap{\, .}
  }
\end{align*}
\normalsize
 \end{proof}

 \begin{remark}
   In our search, we have only optimized over the number of nodes in the tree, not over the number of facets that are involved.
   This accounts for the different numbers of relevant facets in Theorems~\ref{thm:1039}--\ref{thm:3513}.
   We have not attempted to find the minimal number of relevant facets in each case.
  \end{remark}

 \begin{remark}
   For none of the four prismatoids have we found non-realizability certificates that only use the ``first'' and ``last'' layer~\cite{prismatoids2019}.
 \end{remark}

 \begin{remark}
  We would like to insist on the fact that irrespectively of the machine computations that have gone into \emph{finding} these certificates,
  the fact that they are \emph{valid} is perfectly checkable by humans.
\end{remark}

\subsection{Jockusch's $3$-spheres}

Recall Jockusch's family of $3$-spheres, whose member $\Delta^3_6$ was shown to be non-realizable in Proposition~\ref{prop:jockusch}.

\begin{theorem}
  \label{thm:jockusch}
 For $n\ge5$, no member $\Delta^3_n$ of Jockusch's family of $3$-spheres is polytopal.
\end{theorem}

\begin{proof}
  The facets of~$\Delta^3_n$ are listed explicitly in \cite[Lemma~3.1]{novik2020new},
  and among them we find a certain combinatorial $3$-ball~$\pm B_n^{3,1}$.
  Moreover, this $3$-ball is the only subcomplex of $\Delta^3_n$ that gets deleted in the inductive construction of~$\Delta^3_{n+1}$ from~$\Delta^3_n$.
  Therefore, to prove our claim for $n\ge6$ it suffices to check that
  no determinant in the Pl\"ucker certificate from Proposition~\ref{prop:jockusch} involves the facets of $\pm B_6^{3,1}$,
  which are
  \begin{center}\ttfamily \small
    $[ 1 2 5 6 ]$\  $[ \bar{1} \bar{2} \bar{5} \bar{6} ]$\  $[ \bar{1} \bar{2} 5 6 ]$\  $[ 1 2 \bar{5} \bar{6} ]$\  $[ 2 3 5 6 ]$\  $[ \bar{2} \bar{3} \bar{5} \bar{6} ]$\  $[ \bar{2} \bar{3} 5 6 ]$\  $[ 2 3 \bar{5} \bar{6} ]$\  $[ 3 4 5 6 ]$ 
    
 $[ \bar{3} \bar{4} \bar{5} \bar{6} ]$\     $[ \bar{3} \bar{4} 5 6 ]$\  $[ 3 4 \bar{5} \bar{6} ]$\  $[ 1 \bar{4} 5 6 ]$\  $[ \bar{1} 4 \bar{5} \bar{6} ]$\  $[ 1 \bar{4} \bar{5} 6 ]$\  $[ \bar{1} 4 5 \bar{6} ]$\  $[ 1 \bar{4} \bar{5} \bar{6} ]$\  $[ \bar{1} 4 5 6 ]$\rlap{\ .}
  \end{center}
  
  Unfortunately, the certificate of Proposition~\ref{prop:jockusch} \emph{does} contain, for example, the forbidden facet $[1\bar4\bar5\bar6]$!
  We therefore need a different certificate. The provably smallest useful tree certificate is
  \small
  \begin{align*}
\mathtt{
[ 1 2 \bar{6} \bar{4} \bar{2} ] ([ \bar{1} \bar{2} \bar{6} \bar{4} 1 ] ([ 2 3 \bar{6} \bar{4} 5 ] ([ \bar{1} 2 4 \bar{5} 6 ] ([ 2 3 \bar{6} \bar{4} 5 ] ([ \bar{1} 2 \bar{3} 4 3 ] ([ 1 2 \bar{6} \bar{4} \bar{1} ] ([ 1 2 4 \bar{3} \bar{6} ] \,\bblue{\Gamma( 1 2 4 | 3 \bar{3} \bar{4} 5 )} 
}\ \: & \\ \mathtt{
 {}+ [ 1 2 4 \bar{3} 3 ] \bblue{\,\Gamma( 1 2 4 | \bar{3} \bar{4} 5 \bar{6} )}) 
} & \\ \mathtt{
 {}+ [ 1 2 4 \bar{3} 3 ]\,[ 1 2 4 \bar{3} \bar{4} ] \bblue{(-\Gamma( 1 2 \bar{6} | \bar{1} 4 \bar{4} 5 ))}) 
} & \\ \mathtt{
 {}+ [ 1 2 4 \bar{3} 3 ]\,[ 1 2 4 \bar{3} \bar{4} ]\,[ 1 2 \bar{6} \bar{4} 5 ] \bblue{\,\Gamma( \bar{1} 2 4 | 1 3 \bar{3} \bar{6} )}) 
} & \\ \mathtt{
 {}+ [ 1 2 4 \bar{3} 3 ]\,[ 1 2 4 \bar{3} \bar{4} ]\,[ 1 2 \bar{6} \bar{4} 4 ]\,[ \bar{1} 2 \bar{3} 4 3 ] \bblue{\,\Gamma( 2 5 \bar{6} | 1 \bar{1} 3 \bar{4} )}) 
} & \\ \mathtt{
 {}+ [ 1 2 4 \bar{3} 3 ]\,[ 1 2 4 \bar{3} \bar{4} ]\,[ 1 2 \bar{6} \bar{4} 5 ]\,[ \bar{1} 2 \bar{3} 4 \bar{6} ]\,[ 2 3 \bar{6} \bar{4} 5 ] \bblue{\,\Gamma( \bar{1} 2 4 | 1 3 \bar{5} 6 )}) 
} & \\ \mathtt{
 {}+ [ \bar{1} 2 \bar{3} 4 1 ]\,[ 1 2 4 \bar{3} 3 ]\,[ 1 2 4 \bar{3} \bar{4} ]\,[ 1 2 \bar{6} \bar{4} 5 ]\,[ \bar{1} 2 4 \bar{5} 6 ]\,[ 2 3 \bar{6} \bar{4} 5 ] \bblue{(-\Gamma( 2 3 \bar{6} | \bar{1} 4 \bar{4} 5 ))}) 
} & \\ \mathtt{
 {}+ [ 1 2 4 \bar{3} 3 ]\,[ 1 2 3 5 \bar{6} ]\,[ 1 2 4 \bar{3} \bar{4} ]\,[ 1 2 \bar{6} \bar{4} 4 ]\,[ \bar{1} 2 \bar{3} 4 3 ]\,[ \bar{1} 2 4 \bar{5} 6 ]\,[ 2 3 \bar{6} \bar{4} 5 ] \bblue{\,\Gamma( \bar{1} \bar{4} \bar{6} | 1 2 \bar{2} 5 )}) 
} & \\ \mathtt{
 {}+ [ 1 2 4 \bar{3} 3 ]\,[ 1 2 3 5 \bar{6} ]\,[ 1 2 4 \bar{3} \bar{4} ]\,[ 1 2 \bar{6} \bar{4} 4 ]\,[ \bar{1} \bar{2} \bar{6} \bar{4} 2 ]\,[ \bar{1} 2 \bar{3} 4 3 ]\,[ \bar{1} 2 4 \bar{5} 6 ]\,[ 2 3 \bar{6} \bar{4} 5 ] \bblue{(-\Gamma( 1 \bar{4} \bar{6} | \bar{1} 2 \bar{2} 5 ))}
    }
    \rlap{\ ,}    
  \end{align*}
\normalsize
which simplifies to
\small
\begin{align*}
  \mathtt{
  [ \bar{1} 2 \bar{3} 4 1 ]\,[ \bar{1} \bar{2} \bar{6} \bar{4} 1 ]\,[ 1 2 \bar{6} \bar{4} \bar{2} ]\,[ 1 2 4 \bar{3} 3 ]\,[ 1 2 4 \bar{3} \bar{4} ]\,[ 1 2 \bar{6} \bar{4} 5 ]\,[ 2 3 \bar{6} \bar{4} \bar{1} ]\,[ \bar{1} 2 4 \bar{5} 6 ]\,[ 3 4 \bar{6} 5 2 ]\,[ 2 3 \bar{6} \bar{4} 5 ]
          } &
    \\
    \mathtt{{}
  + [ \bar{1} 2 \bar{3} 4 1 ]\,[ \bar{1} \bar{2} \bar{6} \bar{4} 1 ]\,[ 1 2 \bar{6} \bar{4} \bar{2} ]\,[ 1 2 4 \bar{3} 3 ]\,[ 1 2 4 \bar{3} \bar{4} ]\,[ 1 2 \bar{6} \bar{4} 5 ]\,[ \bar{1} 3 \bar{6} 5 2 ]\,[ \bar{1} 2 4 \bar{5} 6 ]\,[ 2 3 \bar{6} \bar{4} 4 ]\,[ 2 3 \bar{6} \bar{4} 5 ]
          } &
    \\
    \mathtt{{}
  + [ \bar{1} 2 4 \bar{5} 1 ]\,[ \bar{1} \bar{2} \bar{6} \bar{4} 1 ]\,[ 1 2 \bar{6} \bar{4} \bar{2} ]\,[ 1 2 4 \bar{3} 3 ]\,[ 1 2 4 \bar{3} \bar{4} ]\,[ 1 2 \bar{6} \bar{4} 5 ]\,[ 2 3 6 4 \bar{1} ]\,[ \bar{1} 2 \bar{3} 4 \bar{6} ]\,[ 2 3 \bar{6} \bar{4} 5 ]^2 \!\!
          } &
    \\
    \mathtt{{}
  + [ 1 2 6 4 \bar{1} ]\,[ \bar{1} \bar{2} \bar{6} \bar{4} 1 ]\,[ 1 2 \bar{6} \bar{4} \bar{2} ]\,[ 1 2 4 \bar{3} 3 ]\,[ 1 2 4 \bar{3} \bar{4} ]\,[ 1 2 \bar{6} \bar{4} 5 ]\,[ 2 3 4 \bar{5} \bar{1} ]\,[ \bar{1} 2 \bar{3} 4 \bar{6} ]\,[ 2 3 \bar{6} \bar{4} 5 ]^2 \!\!
          } &
    \\
    \mathtt{{}
  + [ 1 2 \bar{6} \bar{4} \bar{1} ]\,[ \bar{1} \bar{2} \bar{6} \bar{4} 1 ]\,[ 1 2 \bar{6} \bar{4} \bar{2} ]\,[ 1 2 4 \bar{3} 3 ]\,[ 1 2 4 \bar{3} 5 ]\,[ 1 2 \bar{6} \bar{4} 4 ]\,[ \bar{1} 2 \bar{3} 4 3 ]\,[ \bar{1} 2 4 \bar{5} 6 ]\,[ 2 3 \bar{6} \bar{4} 5 ]^2 \!\!
  } &
    \\
    \mathtt{{}
  + [ 1 2 \bar{6} \bar{4} \bar{1} ]\,[ \bar{1} \bar{2} \bar{6} \bar{4} 1 ]\,[ 1 2 \bar{6} \bar{4} \bar{2} ]\,[ 1 2 \bar{4} 3 4 ]\,[ 1 2 4 \bar{3} 5 ]\,[ 1 2 4 \bar{3} \bar{6} ]\,[ \bar{1} 2 \bar{3} 4 3 ]\,[ \bar{1} 2 4 \bar{5} 6 ]\,[ 2 3 \bar{6} \bar{4} 5 ]^2 \!\!
          } &
    \\
    \mathtt{{}
  + [ 1 2 \bar{6} \bar{4} \bar{1} ]\,[ \bar{1} \bar{2} \bar{6} \bar{4} 1 ]\,[ 1 2 \bar{6} \bar{4} \bar{2} ]\,[ 1 2 3 5 4 ]\,[ 1 2 4 \bar{3} \bar{4} ]\,[ 1 2 4 \bar{3} \bar{6} ]\,[ \bar{1} 2 \bar{3} 4 3 ]\,[ \bar{1} 2 4 \bar{5} 6 ]\,[ 2 3 \bar{6} \bar{4} 5 ]^2 \!\!
          } &
    \\
    \mathtt{{}
  + [ 1 2 \bar{6} \bar{4} \bar{1} ]\,[ 1 2 \bar{6} \bar{4} \bar{2} ]\,[ 1 2 4 \bar{3} 3 ]\,[ 1 2 3 5 \bar{6} ]\,[ 1 2 4 \bar{3} \bar{4} ]\,[ 1 2 \bar{6} \bar{4} 4 ]\,[ \bar{1} 2 \bar{3} 4 3 ]\,[ \bar{1} 2 4 \bar{5} 6 ]\,[ \bar{1} \bar{2} \bar{6} \bar{4} 5 ]\,[ 2 3 \bar{6} \bar{4} 5 ]
          } &
    \\
    \mathtt{{}
  + [ 1 2 \bar{6} \bar{4} \bar{1} ]\,[ 1 2 4 \bar{3} 3 ]\,[ 1 2 3 5 \bar{6} ]\,[ 1 2 4 \bar{3} \bar{4} ]\,[ 1 2 \bar{6} \bar{4} 4 ]\,[ 1 \bar{2} \bar{4} 5 \bar{6} ]\,[ \bar{1} \bar{2} \bar{6} \bar{4} 2 ]\,[ \bar{1} 2 \bar{3} 4 3 ]\,[ \bar{1} 2 4 \bar{5} 6 ]\,[ 2 3 \bar{6} \bar{4} 5 ]
          } &
    \\
    \mathtt{{}
  + [ \bar{1} \bar{2} \bar{6} \bar{4} 1 ]\,[ 1 2 \bar{6} \bar{4} \bar{2} ]\,[ 1 2 4 \bar{3} 3 ]\,[ 1 2 4 \bar{3} \bar{4} ]\,[ 1 2 \bar{6} \bar{4} 4 ]\,[ 1 2 \bar{6} \bar{4} 5 ]\,[ \bar{1} 2 \bar{3} 4 3 ]\,[ \bar{1} 3 \bar{6} 5 2 ]\,[ \bar{1} 2 4 \bar{5} 6 ]\,[ 2 3 \bar{6} \bar{4} 5 ]
          } &
    \\
    \mathtt{{}
  + [ \bar{1} \bar{2} \bar{6} \bar{4} 1 ]\,[ 1 2 4 \bar{3} 3 ]\,[ 1 2 3 5 \bar{6} ]\,[ 1 2 4 \bar{3} \bar{4} ]\,[ 1 2 \bar{6} \bar{4} 4 ]\,[ 1 2 \bar{6} \bar{4} 5 ]\,[ \bar{1} \bar{2} \bar{6} \bar{4} 2 ]\,[ \bar{1} 2 \bar{3} 4 3 ]\,[ \bar{1} 2 4 \bar{5} 6 ]\,[ 2 3 \bar{6} \bar{4} 5 ]
   }
    \rlap{\ .}
\end{align*}
\normalsize
None of these terms contains any forbidden facet,
and the minimality of this certificate is guaranteed by the fact that the integer linear program of Section~\ref{sec:ip} was solved to optimality.
The non-realizability of~$\Delta^3_5$ must be checked separately.
We omit its certificate, which needs six Pl\"ucker polynomials.
\end{proof}

\begin{remark}
  Jockusch's non-realizable $3$-sphere $\Delta^3_5$ with $10$~vertices has the same $f$-vector $(10, 40, 60, 30)$ as the boundary complex of the $4$-polytope
  $P=\conv\big(\Diamond^4\cup\{\pm\boldsymbol{1}\}\big)$
  obtained from the cross-polytope~$\Diamond^4$, but these $3$-spheres are not combinatorially isomorphic.
  The facets of the non-realizable $\Delta^3_5$ are
\begin{quote} \small %\color{blue}
  $[1 2 3 5]$ $[1 2 3 \bar{4}]$ $[1 2 4 5]$ $[1 2 \bar{3} 4]$ $[1 2 \bar{3} \bar{5}]$ $[1 2 \bar{4} \bar{5}]$ $[1 \bar{2} 3 5]$ $[1 \bar{2} 3 \bar{4}]$ $[1 \bar{2} \bar{4} 5]$ $[1 \bar{3} 4 5]$ $[1 \bar{3} \bar{4} 5]$ $[1 \bar{3} \bar{4} \bar{5}]$ $[2 3 4 5]$ $[2 3 4 \bar{5}]$ $[2 3 \bar{4} \bar{5}]$ $[\bar{1} 2 4 \bar{5}]$ $[\bar{1} 2 \bar{3} 4]$ $[\bar{1} 2 \bar{3} \bar{5}]$ $[\bar{1} 3 4 5]$ $[\bar{1} 3 4 \bar{5}]$ $[\bar{1} 3 \bar{4} \bar{5}]$ $[\bar{1} \bar{2} 3 5]$ $[\bar{1} \bar{2} 3 \bar{4}]$ $[\bar{1} \bar{2} 4 5]$ $[\bar{1} \bar{2} \bar{3} 4]$ $[\bar{1} \bar{2} \bar{3} \bar{5}]$ $[\bar{1} \bar{2} \bar{4} \bar{5}]$ $[\bar{2} \bar{3} 4 5]$ $[\bar{2} \bar{3} \bar{4} 5]$ $[\bar{2} \bar{3} \bar{4} \bar{5}]$,
\end{quote}
while the facets of the realizable 3-sphere $\partial P$ are
\begin{quote}\small %\color{blue}
$[1 2 \bar{3} 5]$ $[1 2 \bar{3} \bar{4}]$ $[1 2 \bar{4} 5]$ $[1 3 \bar{4} 5]$ $[1 \bar{2} 3 5]$ $[1 \bar{2} 3 \bar{4}]$ $[1 \bar{2} 4 5]$ $[1 \bar{2} \bar{3} 4]$ $[1 \bar{2} \bar{3} \bar{5}]$ $[1 \bar{2} \bar{4} \bar{5}]$ $[1 \bar{3} 4 5]$ $[1 \bar{3} \bar{4} \bar{5}]$ $[2 3 \bar{4} 5]$ $[2 \bar{3} 4 5]$ $[2 \bar{3} \bar{4} \bar{5}]$ $[\bar{1} 2 3 5]$ $[\bar{1} 2 3 \bar{4}]$ $[\bar{1} 2 4 5]$ $[\bar{1} 2 \bar{3} 4]$ $[\bar{1} 2 \bar{3} \bar{5}]$ $[\bar{1} 2 \bar{4} \bar{5}]$ $[\bar{1} 3 4 5]$ $[\bar{1} 3 \bar{4} \bar{5}]$ $[\bar{1} \bar{2} 3 4]$ $[\bar{1} \bar{2} 3 \bar{5}]$ $[\bar{1} \bar{2} 4 \bar{5}]$ $[\bar{1} \bar{3} 4 \bar{5}]$ $[\bar{2} 3 4 5]$ $[\bar{2} 3 \bar{4} \bar{5}]$ $[\bar{2} \bar{3} 4 \bar{5}]$.
\end{quote}

\end{remark}

\subsection{Novik and Zheng's centrally symmetric neighborly $d$-spheres}

In \cite{novik-zheng20}, Novik and Zheng give several constructions of centrally symmetric, highly neighborly $d$-spheres.
They are based on a family~$\Delta^d_n$ of cs-$\lceil\frac d 2\rceil$-neighborly combinatorial $d$-spheres on $2n\ge 2d+2$ vertices,
which arise as the case $i=\lceil\frac d 2\rceil$ of an inductively constructed family~$\Delta_n^{d,i}$ of cs-$i$-neighborly combinatorial $d$-spheres.
Each of \emph{those} contains a certain combinatorial $d$-ball~$B^{d,i-1}_n$,
which is the only part that gets deleted in a step of the inductive construction.
For $d=3$, Novik and Zheng's family $\{\Delta^3_n: n\ge 4\}$ is precisely Jockusch's family from~\cite{jockusch95},
and $B_n^{3,1}$ is precisely the ball~$B^3_n$ mentioned in Theorem~\ref{thm:jockusch}.

\begin{theorem}
  \label{thm:nz4}
  For $n\ge6$, no member $\Delta^4_n$ of Novik and Zheng's family is realizable.
\end{theorem}

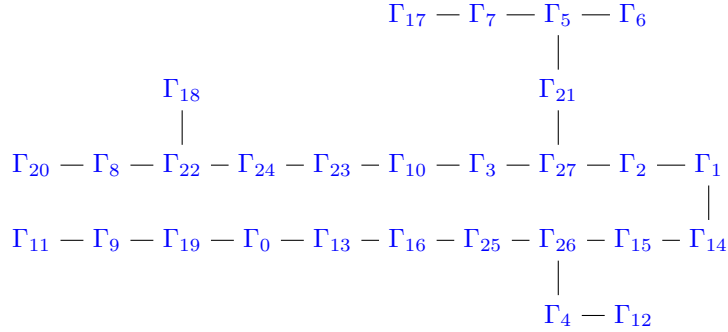
\begin{figure}[htbp]
  \centering
  \begin{tikzpicture}[
    every node/.style={blue}
    ]
    \node (G20) at (0,1) { $\Gamma_{20}$ };
    \node (G8) at (1,1) { $\Gamma_{8}$ };
    \node (G22) at (2,1) { $\Gamma_{22}$ };
    \node (G24) at (3,1) { $\Gamma_{24}$ };
    \node (G23) at (4,1) { $\Gamma_{23}$ };
    \node (G10) at (5,1) { $\Gamma_{10}$ };
    \node (G3) at (6,1) { $\Gamma_{3}$ };
    \node (G27) at (7,1) { $\Gamma_{27}$ };
    \node (G2) at (8,1) { $\Gamma_{2}$ };
    \node (G1) at (9,1) { $\Gamma_{1}$ };

    \node (G11) at (0,0) { $\Gamma_{11}$ };
    \node (G9) at (1,0) { $\Gamma_{9}$ };
    \node (G19) at (2,0) { $\Gamma_{19}$ };
    \node (G0) at (3,0) { $\Gamma_{0}$ };
    \node (G13) at (4,0) { $\Gamma_{13}$ };
    \node (G16) at (5,0) { $\Gamma_{16}$ };
    \node (G25) at (6,0) { $\Gamma_{25}$ };
    \node (G26) at (7,0) { $\Gamma_{26}$ };
    \node (G15) at (8,0) { $\Gamma_{15}$ };
    \node (G14) at (9,0) { $\Gamma_{14}$ };

    \node (G18) at (2,2) { $\Gamma_{18}$ };

    \node (G21) at (7,2) { $\Gamma_{21}$ };
    \node (G5) at (7,3) { $\Gamma_{5}$ };
    \node (G6) at (8,3) { $\Gamma_{6}$ };
    \node (G7) at (6,3) { $\Gamma_{7}$ };
    \node (G17) at (5,3) { $\Gamma_{17}$ };

    \node (G4) at (7,-1) { $\Gamma_{4}$ };
    \node (G12) at (8,-1) { $\Gamma_{12}$ };
    
    \draw (G20)--(G8)--(G22)--(G24)--(G23)--(G10)--(G3)--(G27)--(G2)--(G1)
    --(G14)--(G15)--(G26)--(G25)--(G16)--(G13)--(G0)--(G19)--(G9)--(G11);
    \draw (G22)--(G18);
    \draw (G27)--(G21)--(G5)--(G7)--(G17);
    \draw (G5)--(G6);
    \draw (G26)--(G4)--(G12);
    
  \end{tikzpicture}
  \caption{The minimal Pl\"ucker tree proving the non-realizability of $\Delta^4_7$}
  \label{fig:nz47}
\end{figure}

\begin{proof}
  The construction of~$\Delta^4_n$ in~\cite[Section~3]{novik-zheng20} can be made explicit as follows~\cite{novik-personal}.

  \begin{enumerate}[(1)]
  \item The facets of $B^{4,1}_n$ are

    \begin{enumerate}[(i)]
    \item $\{i,i+1,n-2,n-1,n\}$ and $\{-i,-i-1,n-2,n-1,n\}$ for $1\leq i \leq
      n-4$, and

    \item $\{1, -n+3, n-2, n-1, n\}$, $\{1, -n+3, -n+2, n-1, n\}$, $\{1, -n+3,-n+2,-n+1, n\}$, $\{1, -n+3,-n+2,-n+1, -n\}$,

    \end{enumerate}

  \item The remaining facets of $\Delta^4_n$ are
    \begin{enumerate}[(i)]
    \item $\{i,i+1,\ell-3,\ell-2,\ell\}$,
      $\{i,i+1,\ell-3,\ell-1,\ell\}$,
      $\{-i,-i-1,\ell-3,\ell-2,\ell\}$,

      $\{-i,-i-1,\ell-3,\ell-1,\ell\}$, where
      $1\leq i \leq \ell-5\leq n-5$ (equivalently,
      $6\leq i+5\leq \ell\leq n$);

    \item 
      $\{1, -\ell+4,\ell-3,\ell-2,\ell\}$, $\{1, -\ell+4,\ell-3,\ell-1,\ell\}$,
      $\{1, -\ell+4,-\ell+3,\ell-2,\ell\}$,

      $\{1, -\ell+4, -\ell+2,\ell-1,\ell\}$,
      $\{1, -\ell+4, -\ell+2,-\ell+1,\ell\}$, $\{1, -\ell+4,-\ell+3, -\ell+1,\ell\}$,
      
      $\{-\ell+4, -\ell+3,-\ell+2,\ell-1,\ell\}$, $\{-\ell+4, -\ell+3,-\ell+2, -\ell+1,\ell\}$,
      where $6\leq\ell\leq n$;
    \end{enumerate}

  \item 
    together with the following 10 facets:

    \noindent
$\{-1,2,-3,4,-5\}$,
$\{1,2,-3,4,-5\}$, $\{1,2, 3,4,-5\}$, $\{1,2,3,-4,-5\}$,
$\{1,-2,-3,4,-5\}$,

    \noindent
$\{1,-2, 3,4,-5\}$, $\{1,-2,3,-4,-5\}$,
$\{-1,-2,-3,4,-5\}$, $\{-1,-2, 3,4,-5\}$, $\{-1,-2,3,-4,-5\}$.
\end{enumerate}
The \texttt{polymake} implementation\footnote{available starting from release 4.3} of this construction has successfully passed various consistency checks.

As to the non-realizability of these spheres,
as in Theorem~\ref{thm:jockusch} the case $n=6$ has to be dealt with separately, with an omitted certificate consisting of~$18$ Pl\"ucker polynomials.
There is, in this case, no certificate that avoids the facets of the balls~$\pm B_6^{4,1}$.

For $n=7$, we do find a certificate that avoids the facets of~$\pm B_7^{4,1}$, which are
\begin{quote} \ttfamily\small
  \noindent
  $[ 1 2 5 6 7 ]$\ $[ \bar{1} \bar{2} \bar{5} \bar{6} \bar{7} ]$\ $[ \bar{1} \bar{2} 5 6 7 ]$\ $[ 1 2 \bar{5} \bar{6} \bar{7} ]$\ $[ 2 3 5 6 7 ]$\ $[ \bar{2} \bar{3} \bar{5} \bar{6} \bar
  {7} ]$\ $[ \bar{2} \bar{3} 5 6 7 ]$\ $[ 2 3 \bar{5} \bar{6} \bar{7} ]$\ $[ 3 4 5 6 7 ]$\ $[ \bar{3} \bar{4} \bar{5} \bar{6} \bar{7} ]$\

  \noindent
  $[ \bar{3} \bar{4} 5 6 7 ]$\ $[ 3 4 \bar{5} \bar{6
} \bar{7} ]$\ $[ 1 \bar{4} 5 6 7 ]$\ $[ \bar{1} 4 \bar{5} \bar{6} \bar{7} ]$\ $[ 1 \bar{4} \bar{5} 6 7 ]$\ $[ \bar{1} 4 5 \bar{6} \bar{7} ]$\ $[ 1 \bar{4} \bar{5} \bar{6} 7 ]$\ $[ \bar{1
} 4 5 6 \bar{7} ]$\ $[ 1 \bar{4} \bar{5} \bar{6} \bar{7} ]$\ $[ \bar{1} 4 5 6 7 ]$. 
\end{quote}
It is provably minimal, and uses the $28$ Pl\"ucker polynomials
\small
\begin{align*}
\Gamma_0 &= \bblue{-\Gamma( 1 \bar{2} \bar{3} 6 | 4 \bar{4} 5 \bar{7} )},
&\Gamma_1 &= \bblue{\phantom{+}\Gamma( 1 \bar{2} \bar{3} 6 | \bar{4} 5 \bar{5} \bar{7} )},
&\Gamma_2 &= \bblue{\phantom{+}\Gamma( 1 \bar{2} 5 6 | \bar{3} 4 \bar{4} \bar{7} )},
&\Gamma_3 &= \bblue{-\Gamma( 1 \bar{3} 4 5 | \bar{2} 6 7 \bar{7} )},
\\
\Gamma_4 &= \bblue{-\Gamma( 1 \bar{3} 6 7 | \bar{2} 4 \bar{4} \bar{6} )},
&\Gamma_5 &= \bblue{\phantom{+}\Gamma( \bar{1} \bar{2} \bar{3} 5 | \bar{4} \bar{5} 6 \bar{7} )},
&\Gamma_6 &= \bblue{-\Gamma( \bar{1} \bar{2} \bar{3} \bar{5} | 4 \bar{4} 5 \bar{7} )},
&\Gamma_7 &= \bblue{-\Gamma( \bar{1} \bar{2} \bar{3} 6 | 1 4 5 \bar{5} )},
\\
\Gamma_8 &= \bblue{\phantom{+}\Gamma( \bar{2} 3 5 6 | 1 \bar{3} 4 \bar{7} )},
&\Gamma_9 &= \bblue{-\Gamma( \bar{2} 3 5 6 | \bar{1} \bar{3} 4 \bar{7} )},
&\Gamma_{10} &= \bblue{\phantom{+}\Gamma( \bar{2} \bar{3} 4 5 | 1 \bar{1} \bar{5} \bar{7} )},
&\Gamma_{11} &= \bblue{-\Gamma( \bar{2} \bar{3} 4 6 | 1 \bar{1} 3 5 )},
\\
\Gamma_{12} &= \bblue{\phantom{+}\Gamma( \bar{2} \bar{3} \bar{4} 6 | 1 4 \bar{5} 7 )},
&\Gamma_{13} &= \bblue{-\Gamma( \bar{2} \bar{3} \bar{4} 6 | 1 5 \bar{6} \bar{7} )},
&\Gamma_{14} &= \bblue{\phantom{+}\Gamma( \bar{2} \bar{3} \bar{4} 6 | 1 \bar{5} \bar{6} \bar{7} )},
&\Gamma_{15} &= \bblue{\phantom{+}\Gamma( \bar{2} \bar{3} \bar{4} \bar{6} | 1 \bar{1} 6 7 )},
\\
\Gamma_{16} &= \bblue{-\Gamma( \bar{2} \bar{3} \bar{4} \bar{6} | 1 6 7 \bar{7} )},
&\Gamma_{17} &= \bblue{\phantom{+}\Gamma( \bar{2} \bar{3} 5 6 | 1 \bar{1} 3 \bar{4} )},
&\Gamma_{18} &= \bblue{\phantom{+}\Gamma( \bar{2} \bar{3} 5 6 | 1 3 \bar{4} \bar{5} )},
&\Gamma_{19} &= \bblue{\phantom{+}\Gamma( \bar{2} \bar{3} 5 6 | 1 3 \bar{4} \bar{7} )},
\\
\Gamma_{20} &= \bblue{\phantom{+}\Gamma( \bar{2} \bar{3} 5 6 | \bar{1} 3 4 \bar{4} )},
&\Gamma_{21} &= \bblue{\phantom{+}\Gamma( \bar{2} \bar{3} 5 6 | \bar{1} 4 \bar{4} \bar{7} )},
&\Gamma_{22} &= \bblue{\phantom{+}\Gamma( \bar{2} \bar{3} 5 6 | 3 \bar{4} \bar{5} \bar{7} )},
&\Gamma_{23} &= \bblue{\phantom{+}\Gamma( \bar{2} \bar{3} \bar{5} \bar{7} | 1 4 \bar{4} 5 )},
\\
\Gamma_{24} &= \bblue{-\Gamma( \bar{2} \bar{3} \bar{5} \bar{7} | 1 \bar{4} 5 6 )},
&\Gamma_{25} &= \bblue{-\Gamma( \bar{2} \bar{3} 6 7 | 4 \bar{4} \bar{5} \bar{6} )},
&\Gamma_{26} &= \bblue{\phantom{+}\Gamma( \bar{3} 6 \bar{6} 7 | 1 \bar{2} \bar{4} \bar{5} )},
&\Gamma_{27} &= \bblue{-\Gamma( 4 5 6 \bar{7} | 1 \bar{2} 3 \bar{3} )}
\end{align*}
\normalsize
arranged as in Figure~\ref{fig:nz47}.
Since it contains no facet of $\pm B_7^{4,1}$, it survives the inductive construction.
\end{proof}

\begin{theorem}
  \label{thm:novik-zheng-d}
  \cite{zheng-personal}
  For $n-2\ge d\ge3$, no member $\Delta^d_n$ of Novik and Zheng's family is realizable.
\end{theorem}

\begin{proof}
  By \cite[Proposition~4.1]{novik-zheng20}, each sphere $\Delta^d_n$ occurs as a face link in $\Delta^{d+2}_{n+2}$ and $\Delta^{d+3}_{n+3}$.
  Since links of realizable spheres are realizable, Theorem~\ref{thm:jockusch} yields the proof for all $n\ge d+2$ and $d\ge5$,
  and Theorem~\ref{thm:nz4} directly settles the remaining case $d=4$.
\end{proof}

\begin{remark}
  In \cite{novik2020new}, Novik and Zheng describe several other families of highly neighborly centrally symmetric spheres.
  Once explicit facet descriptions of these are implemented, they can be checked for realizability using the present methods.
\end{remark}

 \section{Acknowledgements}

 It is a pleasure to thank Francisco Santos, Michael Joswig and G\"unter M.~Ziegler for crucial discussions
 and their careful reading and pertinent suggestions.
 Moreover, I am very grateful to Amy Wiebe and Antonio Macchia for finding
 and pointing out an error in a previous version,
 and to Isabella Novik and Heilun Zheng for stimulating discussions, and for pointing out various consequences of
 this work to their families of cs-neighborly spheres.
 
\bibliographystyle{halpha}
\bibliography{pptcert}

\end{document}